\documentclass[11 pt]{amsart}
\usepackage{amsmath,amsfonts,amsthm,mathrsfs,amssymb,bbm}
\usepackage{tikz}
\usepackage{tikz-cd}
\usepackage[margin=1.1in]{geometry}
\usepackage{appendix}

\newtheorem{thm}{Theorem}[section]
\newtheorem{lem}{Lemma}[section]
\newtheorem{prop}{Proposition}[section]
\newtheorem{defn}{Definition}[section]

\usepackage{chngcntr}
\counterwithin{equation}{section}

\newcommand{\bd}{\begin{document}}
\newcommand{\ed}{\end{document}}
\newcommand{\beq}{\begin{equation}}
\newcommand{\eeq}{\end{equation}}
\newcommand{\bp}{\begin{proof}}
\newcommand{\ep}{\end{proof}}
\newcommand{\mytitle}[2]{\begin{center}{\bf {#1}}\\{Geoff Bentsen}\\{#2}\end{center}}

\newcommand{\bq}{\begin{quote}}
\newcommand{\eq}{\end{quote}}
\newcommand{\be}{\begin{enumerate}}
\newcommand{\ee}{\end{enumerate}}
\newcommand{\bi}{\begin{itemize}}
\newcommand{\ei}{\end{itemize}}

\newcommand{\zz}{\mathbb{Z}} 
\newcommand{\rr}{\mathbb{R}} 
\newcommand{\nn}{\mathbb{N}} 
\newcommand{\qq}{\mathbb{Q}} 
\newcommand{\cc}{\mathbb{C}} 
\newcommand{\ff}{\mathbb{F}} 

\usepackage{mathtools} 
\newcommand{\essinf}{\mathrm{ess}\inf}
\newcommand{\dist}{\mathrm{dist}}
\newcommand{\norm}[1]{\Big\|{#1}\Big\|} 
\newcommand{\Ker}{\mathrm{ker}}
\newcommand{\mx}[1]{\begin{smallmatrix} #1 \end{smallmatrix}} 
\newcommand{\supp}{\mathrm{supp}}
\newcommand{\Span}{\mathrm{Span}}
\newcommand{\dil}{\mathrm{Dil}}
\newcommand{\diam}{\mathrm{diam}}
\DeclareMathOperator{\sgn}{sgn}
\newcommand{\C}{\tilde{C}}
\newcommand{\ophi}{\varphi_Q} 
\newcommand{\1}{\mathbbm{1}}

\newcommand{\cA}{\mathcal{A}}
\newcommand{\cB}{\mathcal{B}}
\newcommand{\cC}{\mathcal{C}}
\newcommand{\cE}{\mathcal{E}}
\newcommand{\cI}{\mathcal{I}}
\newcommand{\cL}{\mathcal{L}}
\newcommand{\cM}{\mathcal{M}}
\newcommand{\cR}{\mathcal{R}}
\newcommand{\cT}{\mathcal{T}}
\newcommand{\cS}{\mathcal{S}}
\newcommand{\cZ}{\mathcal{Z}}
\newcommand{\cF}{\mathcal{F}}

\newcommand{\fA}{\mathfrak{A}}
\newcommand{\fF}{\mathfrak{F}} 
\newcommand{\fS}{\mathfrak{S}}
\newcommand{\fw}{\mathfrak{w}}
\newcommand{\fx}{\mathfrak{x}}
\newcommand{\fy}{\mathfrak{y}}
\newcommand{\fz}{\mathfrak{z}}

\DeclareFontFamily{U}{mathx}{\hyphenchar\font45}
\DeclareFontShape{U}{mathx}{m}{n}{
      <5> <6> <7> <8> <9> <10>
      <10.95> <12> <14.4> <17.28> <20.74> <24.88>
      mathx10
      }{}
\DeclareSymbolFont{mathx}{U}{mathx}{m}{n}
\DeclareFontSubstitution{U}{mathx}{m}{n}
\DeclareMathAccent{\widecheck}{0}{mathx}{"71}

\usepackage{scalerel,stackengine}
\stackMath
\newcommand\reallywidehat[1]{%
\savestack{\tmpbox}{\stretchto{%
  \scaleto{%
    \scalerel*[\widthof{\ensuremath{#1}}]{\kern-.6pt\bigwedge\kern-.6pt}%
    {\rule[-\textheight/2]{1ex}{\textheight}}
  }{\textheight}%
}{0.5ex}}%
\stackon[1pt]{#1}{\tmpbox}%
}
\begin{document}

\title{$L^p$ regularity estimates for a class of integral operators with fold blowdown singularities}
\author[G. Bentsen]{Geoffrey Bentsen}
\address{Geoffrey Bentsen \\ Department of Mathematics \\ Northwestern University \\ IL 60208, USA}
\email{geoffrey.bentsen@northwestern.edu}
\thanks{Research supported in part by NSF grant DMS 1764295}
\date{\today}

\begin{abstract}
    We prove sharp $L^p$ regularity results for a class of generalized Radon transforms for families of curves in a three-dimensional manifold associated to a canonical relation with fold and blowdown singularities. The proof relies on decoupling inequalities by Wolff and Bourgain-Demeter for plate decompositions of thin neighborhoods of cones and $L^2$ estimates for related oscillatory integrals.
\end{abstract}
\maketitle

\section{Introduction}

Let $M$ be the family of all lines in $\rr^3$. Given a function $f\in C_0^\infty(\rr^3)$, its X-ray transform is a function defined on $M$ given by
\[
Xf(l)=\int_l f, \qquad l\in M
\]
Since $M$ is a 4 dimensional manifold, recovering $f$ from $Xf$  is an overdetermined problem. It is natural to ask for which 3 dimensional submanifolds $\cF\subset M$ the restriction $X_\cF f=Xf\big|_{\cF}$ can be inverted. We study a class of these restricted X-ray transforms initially formulated in the complex setting by Gelfand and Graev \cite{GeGr68} to give an essentially complete characterization of when inversion is possible.

\begin{defn}[Gelfand Admissibility]\label{Gelfand}
Given a three-dimensional line complex $\cF$, let $\Gamma_Q$ be the cone of lines in $\cF$ through the point $Q$. We say that $\cF$ is {\bf \emph{Gelfand-admissible}} if $\Gamma_P$ is tangent to $\Gamma_Q$ along the line between the points $P$ and $Q$ for every $P$ in the cone $\Gamma_Q$.
\end{defn}

This class of restricted X-ray transforms has been studied by many authors, including Greenleaf and Uhlmann who, in \cite{GrUh89}, showed that Gelfand admissibility, along with the condition that the cone of lines through each point is curved, is sufficient for the inversion of $X_\cF$, extending the results of Gelfand-Graev to the real setting. Various estimates have been proven for this collection of restricted X-ray transforms. For instance, $L^2$-Sobolev estimates were proven by Greenleaf-Uhlmann in \cite{GrUh90B}, and $L^p\to L^q$ estimates were proven by Greenleaf-Seeger in \cite{GrSe94}. In this paper we are interested in finding $L^p$-Sobolev estimates for $X_{\cF}$ and similar operators. It is instructive to look at the following model case. 

Let $I$ be a compact interval and suppose that $\gamma:I\to \rr^2$ is a smooth regular curve with nonvanishing curvature (i.e. $\gamma'(s),\gamma''(s)\ne 0$). For a Schwartz function $f\in\mathcal{S}(\rr^3)$ and $\alpha\in I$ define 
\[
\cA f(x',\alpha)=\int_1^2 f(x'+s\gamma(\alpha),s) \chi_1(s)\chi_2(\alpha) \, ds,
\]

where $\chi_1$ and $\chi_2$ are smooth real-valued functions supported in the interior of $[1,2]$ and $I$ respectively. Pramanik and Seeger, in \cite{PrSe06}, proved that for sufficiently small $p$ the operator $\cA$ maps boundedly from $L^p(\rr^3)$ into $L^p_{1-1/p}(\rr^3)$, where $L^p_s(\rr^3)$ is the standard Sobolev space on $\rr^3$ with respect to Lebesgue measure. This result was proven by studying dyadic decompositions of the adjoint operator $\cA^*$ and using $\ell^p$-decoupling inequalities for the cone, originally proven by Wolff in \cite{Wo00} and later extended to the optimal range by Bourgain-Demeter \cite{BoDe15}. Applying the Bourgain-Demeter decoupling result yields the boundedness of $\cA^*$ from $L^p(\rr^3)$ into $L^p_{1/p}(\rr^3)$ for $p>4$ (and hence the boundedness of $\cA$ from $L^{p'}(\rr^3)$ into $L^{p'}_{1/p}(\rr^3)$ for the same range of $p$, where $p'$ is the dual exponent $p'=1-\tfrac1p$). This estimate is the best possible for the range of $p$, although it is unknown whether the range of $p$ can be extended to include $p=4$. 

 A generalization of the main result of \cite{PrSe06} is suggested by \cite{PrSe19}. In this work, Pramanik and Seeger proved a gain of $1/p$ derivatives in $L^p$ for a class of integral operators in $\rr^3$ with folding canonical relations, generalizing their previous result in \cite{PrSe07}, which considered averages over translations of curves in $\rr^3$. We use similar techniques to Pramanik and Seeger to generalize the results of \cite{PrSe06} to more general integral operators associated to fold and blowdown singularities. This class of integral operators includes the adjoints of generic Gelfand-admissible restricted X-ray transforms, and also subsumes the main result of \cite{Be19} for averaging operators over curves in the Heisenberg group.

To define our class of integral operators, we recall the double fibration formalism of Gelfand and Helgason \cite[p. 4]{He99} (cf. \cite{GrUh89,PrSe19}). Let $\Omega_L,\Omega_R$ be three-dimensional manifolds and consider families of curves $\cM_x\subset \Omega_R$ parametrized by and smoothly depending on $x\in\Omega_L$. Let $d\sigma_x$ be the arclength measure on $\cM_x$, and $\chi\in C_c^\infty(\rr^3\times\rr^3)$. We define the generalized Radon transform operator $\cR: C_c^\infty(\Omega_R)\to C^\infty(\Omega_L)$ by
\[
\cR f(x)=\int_{\cM_x} f(y) \chi(x,y) \, d\sigma_x(y).
\]
We assume that $\cM_x$ are sections of a manifold $\cM\subset \Omega_L\times\Omega_R$, so that the projections 
\begin{equation}
\begin{tikzcd}\label{Mprojection}
& \arrow[dl] \cM \arrow[dr] & \\
\Omega_L & & \Omega_R
\end{tikzcd}
\end{equation}
have surjective differentials; note this ensures that $\cR$ is bounded on $L^1$ and $L^\infty$. The surjectivity assumption on the projections \eqref{Mprojection} also ensures that $\cM_x$ and $\cM_y=\{x\in\Omega_L \ : \ (x,y)\in\cM\}$ are smooth immersed curves in $\Omega_R$ and $\Omega_L$ respectively. 

The operator $\cR$ can be realized as a Fourier integral operator of order $-1/2$ belonging to the H\"ormander class $I^{-\frac12}(\Omega_L,\Omega_R; (N^*\cM)')$, where
\[
(N^*\cM)'=\big\{(x,\xi,y,\eta) \ : \ (x,\xi,y,-\eta)\in N^*\cM\big\}
\]
with $N^*\cM$ the conormal bundle of $\cM$. The assumptions on the projections \eqref{Mprojection} imply that
\[
\cC=(N^*\cM)'\subset(T^*\Omega_L\setminus 0_L)\times (T^*\Omega_R\setminus 0_R),
\]
where $0_L$ and $0_R$ are the zero sections of the cotangent spaces $T^*\Omega_L$ and $T^*\Omega_R$ respectively. Moreover, $\cC$ is a homogeneous canonical relation, i.e. if $\omega_L$ and $\omega_R$ are the canonical two-forms on $T^*\Omega_L$ and $T^*\Omega_R$ respectively, then $\cC$ is Lagrangian with respect to $\omega_L-\omega_R$. As is known from the theory of Fourier integral operators (see \cite{Ho71,Ph95}) the $L^2$-Sobolev regularity properties of $\cR$ are governed by by the geometry of the projections
\begin{equation}\label{canonreldiagram}
\begin{tikzcd}
& \arrow[dl,"\pi_L" '] \cC \arrow[dr,"\pi_R"] & \\
T^*\Omega_L & & T^*\Omega_R
\end{tikzcd}
\end{equation}
This microlocal point of view is due to Guillemin-Sternberg \cite{GuSt77}. Since $\cC$ is Lagrangian the ranks of the differentials $(D\pi_L)_P$ and $(D\pi_R)_P$ are equal; in particular this implies that $(D\pi_L)_P$ is invertible if and only if $(D\pi_R)_P$ is invertible (see \cite{Ho71}). For averaging operators over curves in dimensions larger than 2 the projections $\pi_L$ and $\pi_R$ fail to be diffeomorphisms, meaning that for every point $(x,y)\in\cM$ there is a $P=(x,\xi,y,\eta)\in(N^*\cM)'$ such that $(D\pi_L)_P$ and $(D\pi_R)_P$ are not invertible. However, we can restrict how singular the maps $\pi_L$ and $\pi_R$ are on $\cC$. Following the survey papers \cite{Co97} and \cite{GrSe02}, we recall the definitions of a Whitney fold and a blowdown.

\begin{defn}\label{foldblowdown}
    Suppose $g:X\to Y$ is a $C^\infty$ map between $C^\infty$ manifolds of corank $\le 1$ such that $d(\det(dg)_P\ne 0$ for every $P\in X$ such that $\det(dg)_P=0$. By the implicit function theorem the set $\mathcal{L}=\{P\in X \ : \ \det (dg)_P=0\}$ is thus an immersed hypersurface. 
    We say $V$, a nonzero smooth vector field on $X$, is a {\bf\emph{kernel field}} of $g$ if $V|_P\in \Ker(dg)_P$ for all $P\in \mathcal{L}$.   
    
    We say $g$ is a {\bf\emph{Whitney fold}} if 
    for every kernel field $V$ of $g$ and every $P\in\mathcal{L}$ we have $V(\det dg)\ne 0$ at $P$. 
    
    We say $g$ is a {\bf\emph{blowdown}} if every kernel field $V$ of $g$, when restricted to $\cL$, is everywhere tangential to $\mathcal{L}$. Note this implies that $V^k(\det dg)\big|_P=0$ for all $k\in\nn$ and all $P\in\mathcal{L}$. 
\end{defn}
In \cite{PrSe19}, Pramanik and Seeger proved $\cR$ maps $L^p(\rr^3)$ into $L^p_{1/p}(\rr^3)$ boundedly for $p>4$ for a class of operators where the only singularities on $\pi_L$ and $\pi_R$ are Whitney folds. They conjectured that only the Whitney fold assumption on $\pi_L$ is necessary for their result. In this paper we consider a ``worst'' case, where $\pi_R$ is instead a blowdown. 

\begin{thm}\label{mainthm}
    Let $\cM\subset \Omega_L\times\Omega_R$ be a four-dimensional manifold such that the projections $\cM\to\Omega_L$ and $\cM\to\Omega_R$ are submersions. Assume that the only singularities on $\pi_L:(N^*\cM)'\to T^*\Omega_L$ are Whitney folds, and that $\pi_R:(N^*\cM)'\to T^*\Omega_R$ is a blowdown. Let $\cL$ be the conic submanifold on which $d\pi_L$ and $d\pi_R$ drop rank by one, and let $\varpi$ be the projection of $(N^*\cM)'$ onto the base $\cM$. Suppose that the restriction of $\varpi$ to $\cL$,
    \begin{equation}\label{varpirestriction}
        \varpi:\cL\mapsto \cM
    \end{equation}
    is a submersion. Then $\cR$ extends to a continuous operator
    \[
    \cR \ : \ L^p_{\rm comp}(\Omega_R)\to L^p_{1/p,{\rm loc}}(\Omega_L), \ 4<p<\infty.
    \]
\end{thm} 

Theorem \ref{mainthm} generalizes the results of \cite{Be19} and \cite{PrSe06}, and the sharpness examples in both papers show that the regularity index $s=1/p$ cannot be improved, and that the result fails for $p<4$. Note that the assumption on the projection $\varpi$ ensures a curvature condition on the fibers of $\cL$, first formulated in \cite{GrSe94}, and proven for $\cR$ in \cite{PrSe19}. This curvature ensures that $\ell^p$-decoupling can be applied.

The layout of this paper is as follows. In Section \ref{examplesection} we introduce some example operators for which Theorem \ref{mainthm} applies. In Sections \ref{initialgeom} and \ref{InitialDecompositionSection} we begin the proof of Theorem \ref{mainthm} by relating it to an estimate of oscillatory integrals in Proposition \ref{decomposedfinal}. This is the main estimate of the paper, proven through the interpolation of a decoupling inequality and an $L^2$ estimate in Sections \ref{DecouplingSection} and \ref{L2section} respectively. While the $L^2$ boundedness of $\cR$ has been established by the work of Greenleaf and Seeger in \cite{GrSe94}, these estimates rely on a Strichartz-type argument that does not yield the quantitative estimates that we need to interpolate with the $\ell^p$-decoupling estimates in Section \ref{DecouplingSection}. The work of Comech in \cite{Co99} establishes these quantitative estimates if $\pi_R$ is of finite type but does not cover the case when $\pi_R$ is a blowdown, which is what we prove in Section \ref{L2section} in a general setting. Fortunately it is not necessary to prove the endpoint $L^2$ estimate (see \cite{PhSt97}) in order to interpolate. Finally in Section \ref{PRSsection} we finish the proof of Theorem \ref{mainthm} with a Calder\'on-Zygmund type estimate proven in \cite{PrRoSe11}.

\section{Some Examples}\label{examplesection}
Now we elaborate on some examples to which Theorem \ref{mainthm} applies. The notation in this section is self-contained. 
\subsection{Averages along curves in $\mathbb{H}^1$}

Define the Heisenberg group $\mathbb{H}^1$ to be $\rr^3$ with the group operation 
\[
x\odot y=\big(x_1+y_1,x_2+y_2,x_3+y_3+\tfrac12(x_1y_2-x_2y_1)\big).
\]
Let $\gamma:[0,1]\to\rr^3$ be a smooth regular curve whose tangent vector is nowhere parallel to $(0,0,1)$, so that without loss of generality we can write $\gamma(t)=(t,\gamma_2(t),\gamma_3(t))$. Let $\mu$ be a smooth measure supported on $\gamma([0,1])$, so for $f\in \cS(\rr^3)$ we define
\[
Af(x)=\int_0^1 f(\gamma(t)^{-1}\odot x) \, d\mu(t).
\]
Secco, in \cite{sec99}, developed a group-invariant notion for higher derivatives of $\gamma$ and formulated two conditions which serve as right- and left-invariant analogues of nonvanishing curvature and torsion. These conditions are
\begin{align}\label{heispiL}
    \det\left(\mx{\gamma_2''(t) & \gamma_3''(t) \\ \gamma_2'''(t) & \gamma_3'''(t)}\right)+\tfrac12(\gamma_2''(t))^2&\ne0 \\
\label{heispiR}
    \det\left(\mx{\gamma_2''(t) & \gamma_3''(t) \\ \gamma_2'''(t) & \gamma_3'''(t)}\right)-\tfrac12(\gamma_2''(t))^2&\ne 0.
\end{align}
In \cite{Be19}, the author showed that if \eqref{heispiL} holds for all $t\in\supp(\chi)$ and \eqref{heispiR} does not hold for any $t\in\supp(\chi)$ then $A$ maps boundedly from $L^p_{comp}(\rr^3)$ into $L^p_{1/p}(\rr^3)$ for $p>4$. Under this condition, the operator $A$ is a Fourier integral operator where $\pi_L$ is a fold and $\pi_R$ is a blowdown. An example of a curve satisfying this condition is $\gamma(t)=(t,t^2,\tfrac16t^3)$. 

We next check that $A$ satisfies the final condition of Theorem \ref{mainthm}. The associated incidence manifold $\cM$ is given by
\[
\cM=\big\{(x,y) \ : \ \Phi(x,y)=0\},
\]
where 
\[
\Phi(x,y)=\Big(\mx{\Phi^1(x,y) \\ \Phi^2(x,y)}\Big)=\Big(\mx{x_2-y_2-\gamma_2(x_1-y_1) \\ x_3-y_3-\gamma_3(x_1-y_1)+\tfrac12x_1\gamma_2(x_1-y_1)-\tfrac12x_2(x_1-y_1)}\Big).
\]
The twisted conormal bundle is given by 
\[
(N^*\cM)'=\big\{\big(x,(\tau\cdot \Phi)_x,y,-(\tau\cdot \Phi)_y\big) \ : \ \Phi(x,y)=0\big\}
\]
and $\cL$ is the submanifold of $(N^*\cM)'$ defined by 
\[
(\tau_2,\tau_3)\perp \big(\gamma_2''(x_1-y_1), \gamma_3''(x_1-y_1)-\tfrac12x_1\gamma_2''(x_1-y_1)\big).
\]
The condition that \eqref{heispiL} holds for all $t\in\supp(\chi)$ and \eqref{heispiR} does not hold for any $t\in\supp(\chi)$ implies $\gamma_2''(t)\ne 0$ for all $t\in\supp(\chi)$, hence the restriction of $(N^*\cM)'$ to $\cL$ amounts to a restriction of the $\tau$ variables to a 1-dimensional linear subspace for each $(x,y)\in\cM$. Thus the projection $\varpi$ defined in Theorem \ref{mainthm} is a submersion and we recover the result from \cite{Be19} that $A:L^p_{comp}(\rr^3)\to L^p_{1/p}(\rr^3)$ for $p>4$. 

\subsection{Restricted X-ray transforms in $\rr^3$}
Let $M$ be the space of lines in $\rr^3$, let $\cF\subset M$ be a 3 dimensional line complex such that the cone of lines through each point is curved, and as defined in the introduction let
\[
X_\cF f(l)=\int_l f, \qquad l\in \cF.
\]
We recall the parametrization of the Lagrangian of $X_\cF$ from a survey paper of Phong, \cite{Ph95} (see also \cite{GrUh89}), to verify that $X_\cF$ satisfies the assumptions of Theorem \ref{mainthm}. As shown in \cite{GrUh89}, the maps $\pi_L$ and $\pi_R$ are respectively a blowdown and a Whitney fold, so we only need to verify that the projection $\varpi|_\cL$ is a submersion in this case. We can view $M$ locally as a submanifold of $T\rr^3$, identifying each line $l$ with a point $P$ and a direction $\gamma$. As a consequence $T^* M$ (resp. its subspace $T^*\cF$) can be identified with the restriction of $T^*(T\rr^3)$ to $T M$ (resp. $T\cF$), viewed as functionals on $T(T\rr^3)$. The defining relation for $X_\cF$ is given by 
\[
Z=\{((P,\gamma),Q) \ : \ (P,\gamma)\in \cF, \ Q\in l\}=\{((P,\gamma),Q) \ : \ (Q-P)\wedge \gamma=0\},
\]
and its twisted conormal bundle, using the formalism above, is given by
\[
N^*Z=\{( ((P,\gamma); \gamma\wedge\tau, (Q-P)\wedge \tau ); (Q ; \tau\wedge \gamma))\big|_{T_l \cF} \ : \ (Q-P)\wedge \gamma=0 \}. 
\]

At this point we use Jacobi fields (see \cite[Ch. 5]{Do92} and \cite{Kl82}) to make a more concrete characterization of $T_l\cF$ and $T_l M$. Again, these calculations are based on the methods used in \cite{Ph95,GrUh89}. Fixing $l=(P,\gamma)\in\cF$, let $e_0=\gamma$ and pick $e_1,e_2$ such that $e_0,e_1,e_2$ form an orthonormal basis of vectors on $\rr^3$. With $s$ parametrizing arclength on $l$, the line $l$ can be deformed to another line in $M$ by
\[
P+s\gamma\mapsto P+s\gamma+(a_1s+b_1)e_1+(a_2s+b_2)e_2,
\]
where $a_i,b_i$ are any constants. Thus the Jacobi fields $e_1,se_1,e_2,se_2$ can be viewed as a basis for $T_l M$. Given a Jacobi field $X(s)=(a_1s+b_1)e_1+(a_2s+b_2)e_2$, we can view the deformation above using the identification $l=(P,\gamma)\in T\rr^3$ as
\[
(P,\gamma)\mapsto (P+X(0),\gamma+X').
\]
Thus a tangent vector in $T_lM$ can be identified as a pair $(X(0),X')$ lying in $T^*(T\rr^3)$. 

The Gelfand admissibility condition (Definition \ref{Gelfand}) states that along the line $l$, the normal space to $\cF$ is proportional to a fixed vector. This implies we can pick a unit Jacobi field $X_4(s)$ that is normal to $T_l\cF$ and is proportional to a fixed vector perpendicular to $\gamma$. Choose $e_1$ to be this vector, and choose $a,b\in\rr$ so that $a^2+b^2=1$ and $X_4(s)=(a-sb)e_1$. Recall that $M$ is a symplectic manifold with symplectic form given by 
\[
\omega\Big( \sum_{i=1}^2 (a_is+b_i)e_i,\sum_{i=1}^2 (c_i s+d_i)e_i\Big)=\sum_{i=1}^2 b_ic_i-a_id_i.
\]
Then using $X_4(s)$ we can form a symplectic basis for $T_l M$, given by
\[
X_4(s)=(a-sb)e_1, \qquad X_3(s)=(as+b)e_1, \qquad X_2(s)=e_2, \qquad X_1(s)=se_2
\]
Writing $l=(P,\gamma)$ the basis for $T_l\cF$ is given by $\{(X_i(0),X')\}_{i=1}^3$. Let $\{\Phi_i\}_{i=1}^3$ be the dual basis in $T_l^*\cF$. 
Parametrizing $Q$ by its distance $t$ from $P$, i.e. $Q=P+t\gamma$, we can rewrite $\tau\wedge \gamma=\tau_1e_1+\tau_2 e_2$ and $\tau\wedge(Q-P)=t(\tau_1e_1+\tau_2 e_2)$. Then the twisted conormal bundle is given by
\[
N^*Z=\{(P,\gamma);-(at+b)\tau_1 \Phi_1-\tau_2\Phi_2-t\tau_2\Phi_3; P+t\gamma ; \tau_1e_1+\tau_2e_2) \ : \ \tau_1,\tau_2,t\in\rr\}.
\]
We can parametrize $(P,\gamma)=\sum_{i=1}^3 \alpha_i(X_i(0),X_i')$; thus we can parametrize $N^*Z$ by $t,\tau_1,\tau_2,$ and $\alpha_i(X_i(0),X_i')$,  $i=1,2,3$. Using this formalism we can describe $d\pi_L$. Using the above parametrization we can identify $\pi_L$ with the map 
\[
(\{\alpha_i(X_i(0),X_i')\}_{i=1}^3,t,\tau_1,\tau_2)\mapsto \Big(\sum_{i=3} \alpha_i(X_i(0),X_i'), -(at+b)\tau_1 \Phi_1-\tau_2\Phi_2-t\tau_2\Phi_3\Big),
\]
and thus analytically
\[
d\pi_L=\begin{pmatrix}
I & 0 \\
0 & \cB
\end{pmatrix},
\]
where 
\[
\cB=\begin{pmatrix}
-a\tau_1 & -(at+b) & 0 \\
0 & 0 & -1 \\
-\tau_2 & 0 & -t
\end{pmatrix}.
\]

The determinant of this matrix is $\tau_2(at+b)$, so if we make the generic assumption that $2at+b\ne 0$, $\cL$ is exactly the subvariety of $N^*Z$ on which $\tau_2=0$.
The projection $\varpi:N^*Z\to Z$ from Theorem \ref{mainthm} maps
\[
( (P,\gamma);-(at+b)\tau_1 \Phi_1-\tau_2\Phi_2-t\tau_2\Phi_3; Q; \tau_1e_1+\tau_2e_2)\mapsto ((P,\gamma) ; Q).
\]
Since $(P,\gamma; Q)$ is parametrized by only $\alpha_i(X_i(0),X_i')$ and $t$ we see that $\varpi|_{\cL}$ is a submersion. Thus again we see that Theorem \ref{mainthm} generalizes the results of \cite{PrSe06} and applies to the adjoints of restricted X-ray transforms for Gelfand-admissible line complexes as long as the cones $\Gamma_Q$ are curved. As Theorem \ref{mainthm} applies to the adjoints of these restricted X-ray transforms, we see that $X_\cF$ maps boundedly from $L^p_{\rm comp}(\Omega_R)\to L^p_{1-1/p,{\rm loc}}(\Omega_L)$ for $1\le p<4/3$, where $\Omega_L,\Omega_R$ are sufficiently small coordinate patches in $M$ and $\rr^3$, respectively.

\section{Initial Setup}\label{initialgeom}

Using basic facts on generalized Radon transforms we can simplify our operator $\cR$. By localization we may assume that the Schwartz kernel of $\cR$ is supported in a small neighborhood of a base point $P^{\circ}=(x^\circ,y^\circ)\in\cM$. On that neighborhood the manifold $\cM$ can be expressed locally by a defining function $\Phi=(\Phi^1,\Phi^2)^\intercal:\Omega_L\times \Omega_R\to\rr^2$. In other words, $\cM=\{(x,y) \ : \ \Phi(x,y)=0\}$ in a neighborhood of $P^\circ$. Thus using the Fourier inversion formula the Schwartz kernel of $\cR$ is given by an oscillatory integral distribution, formally written as
\begin{equation}\label{initialschwartz}
\chi(x,y)\delta\circ\Phi(x,y)=(2\pi)^{-2}\iint e^{i\tau\cdot \Phi(x,y)} \chi(x,y) \, d\tau.
\end{equation}
Following the procedure found in \cite{PrSe19}, by local changes of variables and possible redefinition of $\chi$, we can write $\cR$ locally as the oscillatory integral operator
\[
\cR f(x)=\iint e^{i\tau\cdot (S(x,y_3)-y')}\chi(x,y)f(y)\, d\tau dy.
\]
The twisted conormal bundle associated to $\cR$ is given by
\begin{align*}
(N^*\cM)'=\{(x,\xi,y,\eta) \ : \ y_i=S^i(x,y_3), \ & i=1,2, \ \xi=\tau_1 S^1(x,y_3)+\tau_2 S^2(x,y_3), \\
 & \eta=(\tau_1,\tau_2,-\tau_1 S^1_{y_3}(x,y_3)-\tau_2 S^2(x,y_3))\}.
\end{align*}
Thus parametrizing $(N^*\cM)'$ by the coordinates $(x_1,x_2,x_3,\tau_1,\tau_2,y_3)$, the projection $\pi_L$ mapping  $(N^*\cM)' \to  T^*\Omega_L$ is identified with the map
\[
\tilde{\pi}_L:(x_1,x_2,x_3,\tau_1,\tau_2,y_3)\mapsto (x,\tau_1S^1_x(x,y_3)+\tau_2 S^2_x(x,y_3)).
\]
Then we see
\[
D\tilde{\pi}_L=\begin{pmatrix} I_{3\times 3}  & 0 \\
\partial_{x_i}\partial_{x_j}(\tau\cdot S) & \cB
\end{pmatrix},
\]
where $\cB=\begin{pmatrix} S^1_x, & S^2_x, & (\tau\cdot S)_{xy_3} \end{pmatrix}.$
Thus we see
\[
\det D\tilde{\pi}_L=\det(S_x^1, \, S_x^2, \, \tau_1S_{xy_3}^1+\tau_2 S_{xy_3}^2)=\tau_1\Delta^1+\tau_2\Delta^2,
\]
where
\[
\Delta^i(x,y_3)=\det(S_x^1, \, S_x^2, \, S_{xy_3}^i)\Big|_{x,y_3}, \ i=1,2.
\]
We define $\cL=\{(x,\xi,y,\eta)\in\cC \ : \ \det D\tilde{\pi}_L=0\}$. Then $\cL$ is a conic submanifold of $(N^*\cM)'$ defined by
\[
\tau_1\Delta^1(x,y_3)+\tau_2\Delta^2(x,y_3)=0.
\]
Similarly, we can identify $\pi_R:(N^*\cM)'\to T^*\Omega_R$ with
\[
\tilde{\pi}_R:(x_1,x_2,x_3,\tau_1,\tau_2,y_3)\mapsto \big(S(x,y_3),y_3,\tau,-\big(\tau_1S_{y_3}^1(x,y_3)+\tau_2S_{y_3}^2(x,y_3)\big)\big).
\]
Let $N(x,y_3)=S_x^1(x,y_3)\wedge S_x^2(x,y_3)$. We see that a kernel field for $\tilde{\pi}_R$ is given by
\[
V_R=\langle N(x,y_3), \nabla_x\rangle.
\]
Indeed, we see that $\langle N(x,y_3), (\tau\cdot S_{y_3})_x\rangle=\tau\cdot \Delta$, and thus vanishes on $\cL$. Note this implies that $-\Delta^2 S_{xy_3}^1+\Delta^1 S^2_{xy_3}\in \Span(S_x^1,S_x^2)$. Since $\pi_R$ is a blowdown, $V_R$ is parallel to $\cL$, which implies $V_R^k(\tau\cdot S_{y_3})=0$ on $\cL$ for all $k\ge 1$. 

Next, we will examine the fibers in $T^*\Omega_L$ of $\cL$. Let $\Sigma_x$ be the fibers of $\pi_L(\cL)$, given by
\[
\Sigma_x=\{(\tau\cdot S)_x(x,y_3) \ : \ \tau\cdot \Delta(x,y_3)=0\}=\{\pm\rho \Xi(x,y_3) \ : \ \rho>0\},
\]
where 
\[
\Xi(x,y_3)=-\Delta^2(x,y_3)S_x^1(x,y_3)+\Delta^1(x,y_3)S_x^2(x,y_3).
\]
Then we see two consequences, one related to our assumption on $\varpi$.

\begin{lem}[\cite{PrSe19}, \S \, 3]
    If $\pi_L$ is a fold and $\varpi$ is a submersion, then $|\Delta|\ne 0$ near $\cL$, and $\Sigma_x$ is a two-dimensional cone that has one non-vanishing principal curvature given by 
    \[
    \rho\langle \Xi_{y_3,y_3},N\rangle.
    \]
\end{lem}

\begin{lem}
    The direction normal to $\Sigma_x$ at a point specified by $(y_3,\rho)$ is given by $N(x,y_3)$.
\end{lem}

\bp

Let $a\in\rr^3$ be fixed. The tangent space of $\Sigma_a$ at a point parametrized by $(y_3,\rho)$ is spanned by 
\begin{align*}
    T_1(a,y_3)&=\Xi(a,y_3) \\
    \tilde{T}_2(a,y_3)&=\Xi_{y_3}(a,y_3),
\end{align*}
so a normal vector at a point $(\rho,y_3)$ is given by 
\begin{align*}
    T_1\wedge \tilde{T}_2&=\Xi\wedge\Xi_{y_3} \\
    &=(\Delta^1\Delta^2_{y_3}-\Delta^2\Delta^1_{y_3})(S_x^1\wedge S_x^2) \\
    & \qquad + (\Delta^1 S_x^2-\Delta^2 S_x^1)\wedge(\Delta^1 S_{x y_3}^2-\Delta^2 S^1_{xy_3}).
\end{align*}
Since $-\Delta^2 S_{xy_3}^1+\Delta^1 S_{x y_3}^2\in\Span(S_x^1,S_x^2)$ for fixed $(x,y_3)$, the expression in the final line of the calculation of $T_1\wedge T_2$ is either $0$ or a scalar multiple of the vector $S_x^1\wedge S_x^2=N$, hence the sum is a multiple of $N(a,y_3)$. 

\ep

\section{Initial Decomposition}\label{InitialDecompositionSection}

We localize in $|\tau|$ then localize away from the singular variety $\cL$, following the ideas of Phong and Stein in \cite{PhSt91}. Let $\chi_0\in C_c^\infty(\rr)$ be equal to 1 on $[\tfrac{1}{2},2]$ and supported on $[\tfrac{1}{4},4]$ such that $\sum_{k\in\zz} \chi_0(2^k\cdot)\equiv 1$. For $k\ge 1$ define $\chi_k(|\tau|)=\chi_0(2^{-k}|\tau|)$. For $\varepsilon>0$ and  $0\le \ell\le \ell_0=\lfloor \tfrac{k}{2+\varepsilon}\rfloor$ let
\[
    a_{k,\ell,\pm}(x,y_3,\tau)=\left\{
    \begin{array}{cc}
         &  \chi_0(2^{\ell-k}(\pm \tau\cdot\Delta(x,y_3))) \qquad \ell<\ell_0 \\
         & 1-\sum_{\ell<\ell_0} \chi_0(2^{\ell-k}(\pm \tau\cdot\Delta(x,y_3))) \qquad \ell=\ell_0
    \end{array}\right.
\]
and define
\begin{align}
     \cR_{k,\ell,\pm}f(x)&=\chi(x_1)\int e^{i\tau\cdot\tilde{\Phi}(x,y)} \chi(y)f(y) \chi_k(|\tau|)a_{k,\ell,\pm}(x,y_3,\tau) \, dy\,d\tau. \label{Akl}
\end{align}
We will suppress the dependence on $\pm$. We prove the following estimate.
\begin{prop}\label{decomposedfinal}
    For $p>4$ and all $\varepsilon>0$ there exists $\varepsilon_0(p)>0$ such that for all $\ell\le \ell_0=\lfloor \tfrac{k}{2+\varepsilon}\rfloor$,
    \[
        \|\cR_{k,\ell}\|_{L^p\to L^p}\le C_p 2^{-(k+\ell\varepsilon_0)/p}.
    \]
\end{prop}
This proposition follows by interpolation with $L^2$ estimates, $L^\infty$ estimates, and a decoupling inequality.  Let $\cI$ be a collection of intervals of length $2^{-\ell}$ with disjoint interiors intersecting a small neighborhood of 0. Then for a function $f:\rr^3\to \rr$ supported in the unit cube and any $I\in\cI$, let $f_{I}(y):=f(y)\1_{I}(y_3)$, so that $f=\sum_{I\in\cI} f_I$ with almost disjoint supports in $y_3$. The necessary $L^2$ estimate is the following.

\begin{prop}\label{L2boundA} 
    Let $\cR_{k,\ell}$ be defined as above. For every $\varepsilon>0$,
    \begin{align}
        \|\cR_{k,\ell}\|_{L^2\to L^2} &\lesssim  2^{(\ell-k)/2+\ell\varepsilon}, \qquad  \ell\le \ell_0=\lfloor \tfrac{k}{2+\varepsilon}\rfloor. \label{L2est}
    \end{align}
    Moreover, by almost disjoint supports of the functions $f_I$,
    \begin{align}
        \norm{\sum_{I\in\cI} \cR_{k,\ell} f_I}_{L^2} &\lesssim  2^{(\ell-k)/2+\ell\varepsilon}\Big(\sum_{I\in\cI}\|f_I\|_{L^2}^2\Big)^{1/2}, \qquad  \ell\le  \ell_0, \label{l2L2est}.
    \end{align}
\end{prop}

Proposition \ref{L2boundA} will be proven in Section \ref{L2section} following methods of almost-orthogonality found in the proof of the Calder\'on-Vaillancourt theorem (see \cite{MuSc13}, \S \, 9.2), originally introduced into this context by Phong and Stein \cite{PhSt91}, Cuccagna \cite{Cu97}, and Comech \cite{Co97}. The main estimate in the proof of Theorem \ref{mainthm} is the decoupling inequality.

\begin{prop}\label{decoupling}
    For every $\varepsilon>0$
    \begin{align*}
        \norm{\sum_{I\in\cI} \cR_{k,\ell} f_I}_{L^p}&\lesssim_\varepsilon 2^{\ell(1/2-1/p+\varepsilon)}\Big(\sum_{I\in\cI} \|\cR_{k,\ell} f_I\|_{L^p}^p\Big)^{1/p}+2^{-10k}\|f\|_{L^p}
    \end{align*}
    for $2\le p\le 6$.
\end{prop}

Following a similar approach to \cite{AnClPrSe2018} and \cite{PrSe19}, we prove Proposition \ref{decoupling} in Section \ref{DecouplingSection} using an inductive argument, at each step combining $l^p$ decoupling with suitable changes of variables.

\begin{proof}[Proof that Propositions \ref{L2boundA} and \ref{decoupling} imply Proposition \ref{decomposedfinal}]

    We begin by proving an $L^\infty$ estimate for $\cR_{k,\ell}$, namely that
    \begin{align}
        \sup_{I\in\cI} \|\cR_{k,\ell}f_I\|_\infty&\lesssim 2^{-\ell} \sup_{I\in\cI} \|f_I\|_\infty \label{linftyLinfty}\\
        \|\cR_{k,\ell} f\|_\infty&\lesssim \|f\|_\infty. \label{Linfty}
    \end{align}

  To see \eqref{linftyLinfty} we estimate the Schwartz kernel of $\cR_{k,\ell}$ (call it $R_{k,\ell}(x,y)$) by integrating by parts in the $\tau$ variables, distinguishing the directions $(\Delta^1,\Delta^2)$ and $(-\Delta^1,\Delta^2)$. This shows that $|R_{k,\ell}(x,y)|\le C_N U_1(x,y) U_2(x,y)$, where 
    \begin{align*}
      U_1(x,y)&=\frac{2^{k-\ell}}{(1+2^{k-\ell}|\Delta^1(y_1-S^1)+\Delta^2(y_2-S^2)|)^N} \\
      U_2(x,y)&=\frac{2^k}{(1+2^k|-\Delta^2(y_1-S^1)+\Delta^2(y_2-S^2)|)^N}
     \end{align*}
    We integrate in $y'$ first, then in $y_3$, which is supported in an interval of length $2^{-\ell}$. To prove \eqref{Linfty} the same argument holds, but we integrate over a larger interval in $y_3$.

    Interpolating \eqref{linftyLinfty} with \eqref{l2L2est} we obtain
    \begin{align}
        \Big(\sum_{I\in\cI} \|\cR_{k,\ell} f_I\|_p^p\Big)^{1/p} &\lesssim_\varepsilon 2^{\ell(3/p-1+\varepsilon)}2^{-k/p}\Big(\sum_{I\in\cI} \|f_I\|_p^p\Big)^{1/p}, \qquad 2\le p\le \infty.
    \end{align}
    Combining this estimate with Proposition \ref{decoupling} we obtain
    \begin{equation}\label{epsilonnegative}
        \|\cR_{k,\ell} f\|_p\lesssim_\varepsilon 2^{\ell(\varepsilon+2/p-1/2)}2^{-k/p}\Big(\sum_{I\in\cI} \|f_I\|_p^p\Big)^{1/p} + 2^{-10k}\|f\|_p, \qquad 2\le p\le 6.
    \end{equation}
    Note that the power of $2^\ell$ in \eqref{epsilonnegative} is negative if $4<p\le 6$ and $\varepsilon$ is sufficiently small. A further interpolation with the $L^\infty$ estimate \eqref{Linfty} yields Proposition \ref{decomposedfinal} for $p>4$. 

\end{proof}


\section{Decoupling}\label{DecouplingSection}

We mirror the structure of the decoupling estimates in \cite{PrSe19}, working out a model case first then reducing the general case to the model case by changes of variables. In the model case, the functions $S^i$ are replaced by $\fS^i$ satisfying simplifying assumptions at the origin. Additionally, the blowdown condition in this model case implies some additional assumptions near the origin.

\subsection{A Model Case}\label{modelcase}
Consider $C^\infty$ maps $(w,z_3)\mapsto \fS^i(w,z_3)$ defined on a neighborhood of $[-r,r]^4$ for some $r\in(0,1)$. For $n\in\nn$ define $M_n>0$ such that
\begin{equation}\label{Mn}
M_n\ge 2+\|\fS^1\|_{C^{n+5}([-r,r]^4)}+\|\fS^2\|_{C^{n+5}([-r,r]^4)},
\end{equation}
where the $C^n$ norm is the supremum of all derivatives orders $0$ to $n$. We assume that for $w\in[-r,r]^3$,
\begin{equation}\label{Smodelx}
(\fS^1,\fS^2,\fS^1_{z_3})\Big|_{(w,0)}=(w_1,w_2,w_3);
\end{equation}
we also assume 
\begin{equation}\label{modelsingularity}
\fS^2_{w,z_3}(0,0)=0,
\end{equation}
and 
\begin{equation}\label{modelfold}
\fS^2_{w_3z_3^2}(0,0)=\kappa_0.
\end{equation}
As the functions $\fS^1,\fS^2$ play the part of $S^1,S^2$ in our model case, we can analyze the geometry of the canonical relation associated to $\fS^1,\fS^2$. Define for $i=1,2$ the functions $\Delta^i_\fS=\det(\fS_w^1, \, \fS_w^2, \, \fS_{wz_3}^i)$. In this model case the singularity surface $\cL_\fS$ is given by the restriction $\mu_1\Delta^1_\fS(w,z_3)+\mu_2\Delta^2_\fS(w,z_3)=0$. We can define the analogue of the right projection $\tilde{\pi}_R: (w,\mu,z_3)\mapsto (\fS(w,z_3),z_3,\mu,-(\mu_1\fS^1_{z_3}(w,z_3)+\mu_2\fS^2_{z_3}(w,z_3)))$, and a kernel field for this map at the point $P$ parametrized by $(w,z_3,\mu)$ is given by
\[
V_R(w,z_3)=\langle \fS_w^1(w,z_3)\wedge \fS_w^2(w,z_3), \nabla_w\rangle.
\]
We assume a blowdown on $\tilde{\pi}_R$, i.e. that $V_R$ is parallel to $\cL_\fS$, implying that \[
V_R^N[\mu_1\Delta^1_\fS+\mu_2\Delta^2_\fS]\Big|_{(w,z_3),\mu\perp \Delta_\fS(w,z_3)}=0
\]
for all $N>0$. Since $\fS^1_w(w,0)=e_1$ and $\fS^2_w(w,0)=e_2$, we see that $V_R(w,0)=\partial_{w_3}$. The above conditions imply that
\begin{align}\label{blowdownS2}
\partial_{w_3}^N \fS^2_{w_3z_3}(w,0)&=0, \ \forall N\ge 1 \\
\partial_{w_3}^N \Delta^2_\fS(w,0)&=0, \ \forall N\ge 1.\label{blowdowndelta2}
\end{align}
Recall that the fibers of the singular manifold $\cL_\fS$ are given for fixed $w$ by
\begin{align*}
\tilde{\Sigma}_w&=\{\mu_1\fS^1_w(w,z_3)+\mu_2\fS^2_w(w,z_3) \ : \ \mu_1\Delta^1_\fS(w,z_3)+\Delta^2_\fS(w,z_3)=0\} \\
&=\{\pm\rho\Xi_\fS(w,z_3) \ : \ \rho>0, \ |z_3|\le r\},
\end{align*}
where $\Xi_\fS(w,z_3)$ is given by $-\fS^1_w(w,z_3)\Delta^2_\fS(w,z_3)+\fS^1_w(w,z_3)\Delta^1_\fS(w,z_3)$. 
Thus $\tilde{\Sigma}_0$ is a cone parametrized by $(\rho,z_3)$ given by 
\[
\{\pm \rho\Xi_\fS(0,z_3) \ : \ \rho>0, \ |z_3|\le r\}=:\Sigma.
\]
Recall from Section \ref{initialgeom} that $\fS^1_w\wedge \fS^2_w(0,b)=:N(b)$ is normal to $\Sigma$ at the point $P$ parametrized by $(\rho',b)$. Thus $T_{P}\Sigma$ has an orthogonal basis given by
\begin{align*}
    T_1(b)&=\Xi_\fS(0,b) \\
    T_2(b)&=T_1(b)\wedge N(b).
\end{align*}
 For $A>1$ and $\delta\ll 1$ let $\Pi_{A,b}(\delta)$ be set of $\xi\in\rr^3$ such that 
 \begin{align*}
     A^{-1}\le |\langle \tfrac{T_1(b)}{|T_1(b)|},\xi\rangle|&\le A \\
    |\langle \tfrac{T_2(b)}{|T_2(b)|},\xi\rangle|&\le A\delta \\
    |\langle \tfrac{N(b)}{|N(b)|},\xi\rangle|&\le A\delta^2.
 \end{align*}
 The sets $\Pi_{A,b}(\delta)$ are unions of $A\times A\delta \times A\delta^2$-boxes with long, middle, and short sides parallel to $T_1(b),$ $T_2(b)$, and $N(b)$ respectively. We will refer to $\Pi_{A,b}(\delta)$ as a plate. Because the cone $\tilde{\Sigma}_0$ is curved we can apply decoupling to the plates $\Pi_{A,b}(\delta)$.
 
 \begin{thm}[\cite{BoDe15}]
 Let $\varepsilon>0$ and $A>1$. There exists a constant $C(\varepsilon,A)$ such that the following holds for $0<\delta_1<\delta_0<1$.
 
 Let $B=\{b_\nu\}_{\nu=1}^M$ be a set of points in an interval $J\subset [-1,1]$ of length $\delta_0$ such that $|b_\nu-b_{\nu'}|\ge \delta_1$ for $b_\nu,b_{\nu'}\in B$, $\nu\ne \nu'$. Let $2\le p\le 6$. Let $f_\nu\in L^p(\rr^3)$ such that the Fourier transform of $f_\nu$ is supported in $\Pi_{A,b}(\delta_1)$. Then
 \[
 \norm{\sum_{\nu} f_\nu}_p\le C(\varepsilon,A)(\delta_0/\delta_1)^{1/2-1/p+\varepsilon}\Big(\sum_\nu\|f_\nu\|_p^p\Big)^{1/p}. 
 \]
 \end{thm}

Let $(w,z_3)\mapsto \alpha(w,z_3)$ be a $C^\infty$ function satisfying for $|(w,z_3)|_\infty<r$,
\begin{align}
M_0^{-1}\le |\alpha(w,z_3)|&\le M_0 \\
|\nabla_w\alpha(w,z_3)|&\le M_0
\end{align}

Let $(w,z,\mu)\mapsto \zeta(w,z,\mu)$ belong to a bounded family of $C^\infty$ functions supported where $|(w,z)|_\infty\le r$ and $1/4 \le |\mu|\le 4$.

Let $\cT_{k,\ell}$ be an operator with Schwartz kernel
\begin{equation}
2^{2k}\int e^{i2^k\langle \mu,\fS(w,z_3)-z'\rangle} \eta\big(2^\ell\alpha(w,z_3)(\mu_1\Delta^1_\fS(w,z_3)+\mu_2\Delta^2_\fS(w,z_3))\big)\zeta(w,z)\eta(|\mu|) \, d\mu.
\end{equation}
The operator $\cT_{k,\ell}$ will play the role of $\cR_{k,\ell}$ after a nonlinear change of variables, while $\alpha(w,z_3)$ is introduced in the localization as a byproduct of those changes of variables.

\begin{prop}\label{modeldecouplingstep}
Let $0<\varepsilon\le 1$, $k\gg 1$, $0\le\ell\le k/2$,
\[
\delta_0\in(2^{-\ell(1-\varepsilon)},2^{-\ell\varepsilon}),
\]
and $\delta_0>\delta_1\ge \max\{2^{-\ell(1-\varepsilon/2)},\delta_0 2^{-\ell\varepsilon/4}\}$. Define $\varepsilon_1=(\delta_1/\delta_0)^2$. Let $J$ be an interval of length $\delta_0$ containing $0$, and $\cI_J$ be a collection of intervals of length $\delta_1$ with disjoint interior and whose interiors all intersect $J$. Let $\sigma\in C_c^\infty(\rr^3)$ be supported $(-1,1)^3$ and define $\sigma_{\ell,\varepsilon_1}(w)=\sigma(2^\ell w_1,2^{\ell}w_2,\varepsilon_1^{-1} w_3)$. Then for $2\le p\le 6$, $g\in L^p(\rr^3)$ with $g_I(y)=g(y)\1_I(y_3)$, and any $N\in\nn$,
\begin{align*}
\norm{\sigma_{\ell,\varepsilon_1}\sum_{I\in \cI_J} \cT_{k,\ell}g_{I}}_{p}&\lesssim_{\varepsilon}(\delta_0/\delta_1)^{1/2-1/p+\varepsilon}\left(\sum_{I\in\cI_J} \norm{\sigma_{\ell,\varepsilon_1}\cT_{k,\ell}g_I}_{p}^p\right)^{1/p} \\
& \qquad \qquad +C(\varepsilon,N)2^{-kN}2^{-2\ell}\varepsilon_1\|g\|_p.
\end{align*}
\end{prop}

The idea here is to show that the Fourier transforms of $\sigma_{\ell,\varepsilon_1}\sum_{I\in \cI_J} \cT_{k,\ell}g_{I}$ are concentrated on the plates $\Pi_{A,b_I}(\delta_1)$ for some $b_I\in I$ and some large enough $A>1$. 

\subsubsection{Derivatives of $\fS$ and $\Delta$}

Some approximations will be helpful to write down. For the rest of Section \ref{modelcase} we omit the subscript dependence on $\fS$. Because of \eqref{Smodelx} we may conclude that for any multiindex $\beta$ of length at least 1,
\begin{align}
\partial_w^\beta \fS^1_w\big|_{(w,0)}&=0 \label{S1wderivs} \\
\partial_w^\beta \fS^2_w\big|_{(w,0)}&=0 \label{S2wderivs} \\
\partial_w^\beta \fS^1_{wz_3}\big|_{(w,0)}&=0. \label{S1zwderivs}
\end{align}

For $w\in[-r,r]^3$, 
\begin{align}\label{delta1approx}
\Delta^1(w,0)&=1 \\
\Delta^2(w,0)&=\fS^2_{w_3z_3}(w,0)\label{delta2approx} \\
\Delta^1_{z_3}(0,0)&=\fS^1_{w_3z_3^2}(0,0)\label{delta1z} \\
\Delta^2_{z_3}(0,0)&=\fS^2_{w_3z_3^2}(0,0)=\kappa_0 \label{delta2z}
\end{align}
and thus
\begin{align}\label{Xiapprox}
\Xi(w,0)&=-\Delta^2(w,0)\fS^1_w(w,0)+\Delta^1(w,0)\fS^2_w(w,0)=e_2-\fS^2_{w_3z_3}(w,0)e_1 \\
\Xi_{w_3^n}(0,0)&=-\fS^2_{w_3^{n+1}z_3}(0,0)e_1=0, \qquad n\ge 1 \label{Xiw3approx} \\
\Xi_{z_3}(0,0)&=-\kappa_0e_1+\fS^1_{w_3z_3^2}(0,0)e_2.\label{Xizapprox}
\end{align}
Using these,
\begin{align}\label{T1bestimate}
    T_1(b)&=\Xi(0,b) \\
    &=\Xi(0,0)+b\Xi_{z_3}(0,0)+O(b^2) \\
    &=-\kappa_0be_1+(1+b\fS^1_{w_3z_3^2}(0,0))e_2+O(b^2)
\end{align}
and
\begin{align}\label{Nbestimate}
    N(b)&=\fS^1_w(0,b)\wedge \fS^2_w(0,b) \\
    &=(e_1+be_3+O(b^2))\wedge(e_2+O(b^2)) \nonumber\\
    &=-be_1+e_3+O(b^2).\nonumber
\end{align}
From these we see that
\begin{align} \label{T2bestimate}
T_2(b)=T_1(b)\wedge N(b)=(1+b\fS^1_{w_3z_3^2}(0,b))e_1+\kappa_0b e_2+ be_3 + O(b^2).
\end{align}

Let $\beta=(\beta_{w_1},\beta_{w_2},\beta_{w_3},\beta_{z_3})$ be a multi-index and let $\partial_{(w,z_3)}^{\beta}$ denote a derivative of order $|\beta|=\beta_{w_1}+\beta_{w_2}+\beta_{w_3}+\beta_{z_3}$ in the variables $w,z_3$. By using the upper bounds $M_n$, trilinearity of determinants, and differentiation rules for products we can estimate
\begin{equation}\label{Deltaderiv}
    |\partial_{(w,z_3)}^{\beta}\Delta^i|\le 3^{|\beta|}M_{|\beta|}^3.
\end{equation}
Similarly, by differentiating products,
\begin{equation}\label{Xideriv}
    |\partial_{(w,z_3)}^\beta\Xi|\le 4^{|\beta|}M_{|\beta|}^4.
\end{equation}

\subsubsection{Plate Localization}

\begin{lem}\label{plateestimatelemma}
    Let $\varepsilon>0$, and $\delta_0,\delta_1,\varepsilon_1$ be as in Proposition \ref{modeldecouplingstep}. Assume that $2^{-\ell}\ll r$, $M_02^{-\ell}\le 2^{-10}$, $\tfrac14\le |\mu|\le 4$, $|w'|\le 2^{-\ell}$, $|w_3|\le \varepsilon_1$, $|b|\le \delta_0$, and $|z_3-b|\le \delta_1$.  

    If 
    \begin{equation}\label{modelfoldlocalization}
    |\mu_1\Delta^1(w,z_3)+\mu_2\Delta^2(w,z_3)|\le M_0 2^{-\ell},
    \end{equation}
    then there exists $A(\varepsilon)>1$ such that
    \[
    \mu_1\fS^1_w(w,z_3)+\mu_2\fS^2_w(w,z_3)\in \Pi_{A(\varepsilon),b}(\delta_1).
    \]
    More specifically, 
    \begin{align}\label{T1plateestimate}
    A(\varepsilon)^{-1}|T_1(b|)\le|\langle T_1(b),\mu_1\fS^1_w(w,z_3)+\mu_2\fS^2_w(w,z_3)\rangle|&\le A(\varepsilon)|T_1(b)|\\
    |\langle T_2(b),\mu_1\fS^1_w(w,z_3)+\mu_2\fS^2_w(w,z_3)\rangle|&\le A(\varepsilon)|T_2(b)|\delta_1. \label{tangentplateestimate} \\
    |\langle N(b), \mu_1 \fS^1_w(w,z_3)+\mu_2 \fS^2_w(w,z_3)\rangle|&\le A(\varepsilon)|N(b)|\delta_1^2. \label{normalplateestimate}
    \end{align}
    Note that the constant $A(\varepsilon)$ does not depend on $\delta_0,\delta_1$.
\end{lem}

\bp

Throughout this proof we use Taylor expansions with appropriate error remainders. Therefore, for any $i=1,2,...$ the function $R_i(w,z_3)$ is $C^\infty$ and uniformly bounded by 1.

The estimate in \eqref{T1plateestimate} is clearly true for some $A>1$ independent of $\varepsilon$. We start with the proof of \eqref{normalplateestimate}.
Let $G=\lceil3\varepsilon^{-1}\rceil$. Employing a Taylor expansion about $(w,z_3)=(0,b)$, and reorganizing terms using that $2^{-\ell}\le \delta_1$, $2^{-\ell}\delta_0\le \delta_1^2$, $\varepsilon_1^2\delta_0^2\le \delta_1^2$,and $\varepsilon_1^G\le \delta_1^2$, we see that
\begin{align}\label{normalexpansion}
\langle N(b), \mu\cdot \fS_w(w,z_3)\rangle&=\sum_{n=0}^G\sum_{|\alpha|=0}^{G-n}\langle N(b), \nabla_w \Big((\partial_{z_3})^n(\partial_w)^\alpha[\mu_1 \fS^1+\mu_2 \fS^2]\Big)(0,b)\rangle \\
&\qquad + M_G\delta_1^2 R_1(w,\mu,z_-3) \nonumber \\
& = \langle N(b), \mu_1 \fS_w^1+\mu_2 \fS_w^2(0,b)\rangle \\
& \qquad + (z_3-b)\langle N(b), \mu_1 \fS^1_{wz_3}+\mu_2 \fS^2_{wz_3}(0,b)\rangle \nonumber \\
& \qquad + \sum_{i=1}^2 w_i\langle N(b), \mu_1 \fS_{w w_i}^1+\mu_2 \fS_{w w_i}^2(0,b)\rangle \nonumber \\
& \qquad +I+II+III+M_G\delta_1^2 R_2(w,\mu,z_3), \nonumber
\end{align}
where
\begin{align*}
I&= \sum_{n=1}^G \frac{w_3^n}{n!} \langle N(b), \mu_1 \fS_{w w_3^n}^1+\mu_2 \fS_{w w_3^n}^2(0,b)\rangle \\
II&=\sum_{n=2}^G \frac{w_3^{n-1}(z_3-b)}{n!}\langle N(b), \mu_1 \fS_{w w_3^{n-1} z_3}^1+\mu_2 \fS_{w w_3^{n-1} z_3}^2(0,b)\rangle \\
III&= \sum_{n=2}^G \sum_{i=1}^2 \frac{w_3^{n-1}w_i}{n!} \langle N(b), \mu_1 \fS_{w w_3^{n-1}w_i}^1+\mu_2 \fS_{w w_3^{n-1}w_i}^2(0,b)\rangle.
\end{align*}
Clearly the first term in \eqref{normalexpansion} vanishes by the definition of $N(b)$ (see \eqref{Nbestimate}). The second term in the expansion is
\begin{align*}
    (z_3-b)\langle N(b), \mu_1 \fS^1_{wz_3}+\mu_2 \fS^2_{wz_3}(0,b)\rangle&=(z_3-b)(\mu_1\Delta^1(0,b)+\mu_2\Delta^2(0,b)).
\end{align*}
    Now, since $|w'|,|z_3-b|\le \delta_1$, applying a Taylor expansion and using trilinearity of determinants, and differentiation of products we get
\begin{align*}
    \mu_1\Delta^1(0,b)+\mu_2\Delta^2(0,b)&=\big(\mu_1\Delta^1(w,z_3)+\mu_2\Delta^2(w,z_3)\big) \\
    &\qquad +\sum_{n=1}^G \frac{w_3^n}{n!}\Big(\mu_1\Delta^1_{w_3^n}(w,z_3)+\mu_2\Delta^2_{w_3^n}(w,z_3)\Big)  \\
    &\qquad +3^G M_G^3\delta_1R_3(w,z_3).
\end{align*}
By \eqref{modelfoldlocalization} the first term is bounded by $M_02^{-\ell}$. For each $1\le n\le G$, from \eqref{blowdowndelta2} and \eqref{delta1approx} we have $\Delta^i_{w_3^{n}}(w,0)=0$ for $i=1,2$, and so by trilinearity of determinants, and differentiation of products, expanding about $z_3=0$ we get
\[
|\Delta^i_{w_3^n}(w,z_3)|\le\Big|\Delta^i_{w_3^n}(w,0)+3^n M_n^3 z_3\Big|\le 3^{n} M_{n}^3\delta_0.
\]
Thus
\[
|\mu_1\Delta^1(0,b)+\mu_2\Delta^2(0,b)|\le M_02^{-\ell}+3^G M_{G}^3\varepsilon_1\delta_0+3^G M_{G}^3\delta_1\le 3^{G+1}M_{G}^3\delta_1,
\]
and the second term in \eqref{normalexpansion} is bounded by $3^{G+1}M_{G}^3\delta_1^2$.

Next we deal with the first order $w'$ derivatives in \eqref{normalexpansion}. We approximate about $z_3=0$. For $i=1,2$, using the estimates \eqref{S1wderivs} and \eqref{S2wderivs}, we get
\begin{align*}
|w_i\langle \fS_w^1(0,b)\wedge \fS_w^2(0,b)&,\mu_1\fS^1_{ww_i}(0,b)+\mu_2\fS^2_{ww_i}(0,b)\rangle|\\
&\le|w_i|\Big[|\langle \fS_w^1(0,0)\wedge \fS_w^2(0,0),\mu_1\fS^1_{ww_i}(0,0) \\
&\qquad +\mu_2\fS^2_{ww_i}(0,0)\rangle|  + 3M_0^3bR_4(0,b)\Big] \\
&\le 2^{-\ell}(0+3 M_0^3\delta_0).
\end{align*}
Note that the condition $\delta_1\ge \max\{M_0^2 2^{20-\ell(1-\varepsilon/2},2^{-\ell\varepsilon/4}\delta_0\}$ from Proposition \ref{modeldecouplingstep} implies that $2^{-\ell}\delta_0\le \delta_1^2$. 

Finally, we estimate $I$, $II$, and $III$. All rely on the blowdown condition at the origin.

First we estimate $I$. For all $n\ge 1$, we expand about the origin to obtain
\begin{align*}
    \langle \fS_w^1(0,b)\wedge \fS_w^2(0,b)&,\mu_1\fS_{ww_3^n}^1(0,b)+\mu_2 \fS_{ww_3^n}^2(0,b)\rangle\\
    &= \langle \fS_w^1(0,0)\wedge \fS_w^2(0,0),\mu_1\fS_{ww_3^n}^1(0,0)+\mu_2 \fS_{ww_3^n}^2(0,0)\rangle \\
    &\qquad +b\Big[\det(\fS^1_{wz_3} \ \fS_w^2 \ \mu_1 \fS_{ww_3^n}^1+\mu_2 \fS_{ww_3^n}^2)\Big|_{(0,0)} \\
    &\qquad +\det(\fS^1_{w} \ \fS_{wz_3}^2 \ \mu_1 \fS_{ww_3^n}^1+\mu_2 \fS_{ww_3^n}^2)\Big|_{(0,0)} \\
    &\qquad +\det(\fS^1_{w} \ \fS_{w}^2 \ \mu_1 \fS_{ww_3^nz_3}^1+\mu_2 \fS_{ww_3^nz_3}^2)\Big|_{(0,0)}\Big] \\
    &\qquad +3^2M_n^3 b^2R_5(0,b)
\end{align*}
Using the estimates \eqref{S1wderivs}, \eqref{S2wderivs}, \eqref{Nbestimate}, \eqref{S1zwderivs}, and \eqref{blowdownS2}, we observe
\begin{align*}
    \langle \fS_w^1(0,0)\wedge \fS_w^2(0,0),\mu_1\fS_{ww_3^n}^1(0,0)+\mu_2 \fS_{ww_3^n}^2(0,0)\rangle&=0 \\
    \det(\fS^1_{wz_3} \ \fS_w^2 \ \mu_1 \fS_{ww_3^n}^1+\mu_2 \fS_{ww_3^n}^2)\Big|_{(0,0)}&=0 \\
    \det(\fS^1_{w} \ \fS_{wz_3}^2 \ \mu_1 \fS_{ww_3^n}^1+\mu_2 \fS_{ww_3^n}^2)\Big|_{(0,0)}&=0 \\
\det(\fS^1_{w} \ \fS_{w}^2 \ \mu_1 \fS_{ww_3^nz_3}^1+\mu_2 \fS_{ww_3^nz_3}^2)\Big|_{(0,0)}&=0.
\end{align*}
This implies 
\[
|I|\le 3^2M_G^3\sum_{n=1}^G \frac{\varepsilon_1^n\delta_0^2}{n!}\le 3^3M_G^3\varepsilon_1\delta_0^2\le 3^3 M_G^3 \delta_1^2.
\]
Next we estimate $II$. For $n\ge 2$, we expand about the origin to obtain
\begin{align*}
    \langle \fS_w^1(0,b)\wedge \fS_w^2(0,b), \mu\cdot \fS_{ww_3^{n-1}z_3}(0,b)\rangle &=\det( \fS_w^1 \ \fS_w^2 \ \mu_1\fS^1_{ww_3^{n-1}z_3}+\mu_2\fS^2_{ww_3^{n-1}z_3})\big|_{(0,0)} \\
    &\qquad +3M_G^3bR_6(0,b).
\end{align*}
Thus the calculation from $I$ the determinant vanishes, and thus
\[
|II|\le 3M_G^3\sum_{n=2}^{G} \frac{\varepsilon_1^{n-1}\delta_1\delta_0}{n!}\le 3 M_G^3 \varepsilon_1\delta_1\delta_0\le 3M_G^3\delta_1^2.
\]
Finally we estimate $III$.  Again using the calculations from $I$, for $n\ge 2$ and $i=1,2$
\begin{align*}
\langle \fS_w^1(0,b)\wedge \fS_w^2(0,b), \mu\cdot \fS_{ww_3^{n-1}w_i}(0,b)\rangle&=\langle \fS_w^1(0,0)\wedge \fS_w^2(0,0), \mu\cdot \fS_{ww_3^{n-1}w_i}(0,0)\rangle \\
&\qquad +3M_G^3bR_7(0,b) \\
&=\mu\cdot \fS_{w_iw_3^n}(0,0)+3M_G^3bR_7(0,b) \\
&=3M_G^3bR_7(0,b).
\end{align*}
This implies that 
\[
|III|\le 3 M_G^3\sum_{n=2}^G \frac{\varepsilon_1^{n-1}2^{-\ell}\delta_0}{n!}\le 3M_G^3 \varepsilon_12^{-\ell}\delta_0\le 3M_G^3\delta_1^2.  
\]
Since $|N(b)|\ge 1/2$ this proves \eqref{normalplateestimate} with any $A(\varepsilon)\ge 3^{\lceil3/\varepsilon\rceil+2}M_{\lceil3/\varepsilon\rceil}^3$.

 Having proven \eqref{normalplateestimate}, we prove \eqref{tangentplateestimate}. Using \eqref{T2bestimate},
define 
\[\
T_2^*(b)=(1+b\fS^1_{w_3z_3^2}(0,0))e_1+\kappa_0be_2+be_3
\]
and note that $|T_2(b)-T_2^*(b)|\le M_0\delta_0^2$. Next, we will approximate $\mu$ by the projection of $\mu_1\Delta^1(w,z_3)+\mu_2\Delta^2(w,z_3)$ onto $\cL_{\fS}$. In particular, let 
\[
\mu^\circ=\pm\tfrac{|\mu|}{|\Delta(w,z_3)|}(-\Delta^2(w,z_3),\Delta^1(w,z_3)),
\]
so that $\mu^\circ_1\Delta^1(w,z_3)+\mu^\circ_2\Delta^2(w,z_3)=0$, $|\mu|=|\mu^\circ|$, and where the sign is picked so that
\[
|\mu-\mu^\circ|\le 2|\mu|M_0 2^{-\ell}.
\]
This is possible since $|\mu_1\Delta^1+\mu_2\Delta^2|\le M_02^{-\ell}$ and $|\Delta(w,z_3)|\ne 0$. Then 
\[
\mu^\circ_1\fS^1_w(w,z_3)+\mu^\circ_2\fS^2_w(w,z_3)=\tfrac{|\mu|}{|\Delta(w,z_3)|}\Xi(w,z_3),
\]
and thus
\[
\Big|\mu_1\fS^1_w(w,z_3)+\mu_2\fS^2_w(w,z_3)-\frac{|\mu|}{|\Delta(w,z_3)|}\Xi(w,z_3)\Big|\le |\mu-\mu^\circ||\fS_w|\le 8M_0^2 2^{-\ell}.
\]
We approximate by a Taylor expansion about the origin, using the fact that $\varepsilon_1\delta_0\le \delta_1$, $|w'|\le 2^{-\ell}\le \delta_1$, $\delta_0^2\le \delta_1$, and $\varepsilon_1^G\le \delta_1^2\le \delta_1$. Reorganizing, we obtain
\begin{align*}
    \langle T_2^*(b),\Xi(w,z_3)\rangle&=\sum_{n=0}^G\sum_{|\alpha|=0}^{G-n}\langle T_2^*(b), (\partial_{z_3})^n(\partial_w)^\alpha\Xi\rangle\Big|_{(0,0)} + 4^G M_G^4 \delta_1 R_8(w,z_3) \\
    &=\langle T^*_2(b),\Xi(0,0)\rangle+z_3\langle T_2^*(b),\Xi_{z_3}(0,0)\rangle \\
    &\qquad +\sum_{n=1}^G \frac{w_3^n}{n!}\langle T^*(b),\Xi_{w_3^n}(0,0)\rangle +  4^GM_G^4\delta_1 R_9(w,z_3).
\end{align*}
Using \eqref{Xiapprox}, \eqref{Xiw3approx}, and \eqref{Xizapprox}
\begin{align*}
    \langle T^*_2(b),\Xi(0,0)\rangle&=\kappa_0b \\
    \langle T^*_2(b),\Xi_{z_3}(0,0)\rangle&=-\kappa_0(1+b(\fS^1_{w_3z_3^2}(0,0)-\fS^1_{w_3z_3^2}(0,b))) \\
    \langle T^*_2(b),\Xi_{w_3^n}(0,0)\rangle&=0, \qquad\qquad \qquad  n\ge 1.
\end{align*}
Thus
\begin{align*}
    |\langle T_2^*(b),\Xi(w,z_3)\rangle| &\le\kappa_0\delta_1 + \kappa_0M_0\delta_0^2+4^GM_G^4\delta_1
\end{align*}
and therefore we can estimate
\begin{align*}
|\langle T_2(b),\mu_1\fS^1_w(w,z_3)+\mu_2\fS^2_w(w,z_3)\rangle|&\le M_0\delta_0^2+8M_0^2 2^{-\ell} +\kappa_0\delta_1+\kappa_0M_0\delta_0^2+4^GM_G^4\delta_1 \\
&\le \kappa_0(1+4^{G+2}M_G^4)\delta_1.
\end{align*}
Thus picking 
\begin{equation}\label{Aepsilon}
    A(\varepsilon)\ge \max\{3^{\lceil3/\varepsilon\rceil+2}M_{\lceil3/\varepsilon\rceil}^3,\kappa_0(1+4^{\lceil3/\varepsilon\rceil+2}M_{\lceil3/\varepsilon\rceil}^4)\}
\end{equation} the Lemma is proven.
\ep

\subsection{Proof of Proposition \ref{modeldecouplingstep}}
Fix an $I\in\cI_J$ and pick $b_I\in I$. Let $m_{A,b_I,\delta_1}$ be be a multiplier equal to 1 on  $\Pi_{2A,b_I}(\delta_1)$ which vanishes on $\Pi_{3A,b_I}(\delta_1)$. Let 
\[
\reallywidehat{P_{k,A,b_I,\delta_1}f}(\xi)=m_{A,b_I,\delta_1}(2^k\xi)\hat{f}(\xi).
\]
Then by Bourgain-Demeter decoupling on the cone,
\[
\norm{\sum_{I} P_{k,A,b_I,\delta_1} \cT_{k,\ell}f_I}_p\le C(\varepsilon,A) (\delta_0/\delta_1)^{1/2-1/p+\varepsilon}\Big(\sum_I \big\|\cT_{k,\ell}f_I\big\|_p^p\Big)^{1/p},
\]
for $2\le p\le 6$. 
The Schwartz kernel of the operator $f\mapsto (\mathrm{I}-P_{k,A,b_I,\delta_1}(\delta_1))Tf$ is given by a sum of kernels $\sum_{n=0}^\infty K_{n,k,\ell}(w,z)$, where
\[
K_{n,k,\ell,b_I}(w,z)=2^{2k}\int\int\int e^{2\pi i\Psi(w,v,z,\mu,\xi)} \sigma_{1}(v,z,\mu)\sigma_{n,2}(\xi) \, dv \, d\xi d\mu,
\]
the phase function $\Psi$ is given by
\[
\Psi(w,v,z,\mu,\xi)=\langle w-v,\xi\rangle+2^k\mu\cdot(\fS(v,z_3)-z')
\]
and the symbols $\sigma_{1},\sigma_{n,2}$ are given by
\begin{align*}
\sigma_{1}(v,z,\mu)&=\sigma_{\ell,\varepsilon_1}(v)\chi_1\big(2^\ell\alpha(v,z_3)\mu\cdot\Delta(v,z_3)\big)\zeta(v,z,\mu), \\
\sigma_{n,2}(\xi)&=\big(1-m_{A(\varepsilon),b_I,\delta_1}(2^k\xi)\big)\chi_n(|\xi|)
\end{align*}
Note that the symbol of $K_{n,k,\ell}$ is supported where $|\xi|\sim 2^n$ for $n\ge 1$ (with obvious modifications for $n=0$), $|\mu|\sim 1$, $|v|+|z|\le r$, and for a priori unbounded $w$. 

We prove the following lemma to reduce to the case when $|\xi|\simeq 2^k$.
\begin{lem}\label{knklxilemma}
There exists a constant $C_1>0$ such that $|k-n|>C_1$ implies that for $N>1$
\begin{equation}\label{knklxibound}
|K_{n,k,\ell,b_I}(w,z)|\le C_{N,\varepsilon} \frac{1}{(1+|w|)^4}2^{-N(k+n)}\1_{[-r,r]}(|z|).
\end{equation}
\end{lem}

If $|n-k|<C_1$ we can apply integration by parts using the fact that $2^{-k}\xi$ is bounded away from the plate $\Pi_{A(\varepsilon),b_I}(\delta_1)$ while $\mu\cdot\fS_w$ lies in $\Pi_{A(\varepsilon),b_I}(\delta_1)$ to obtain lower bounds on $|\Psi_v|$. In particular, we prove the following estimate.
\begin{lem}\label{offplatelemma}
    If $|n-k|<C_1$ then
    \[
        |K_{n,k,\ell,b_I}(w,z)|\le C_{\varepsilon}2^{-11k}\frac{1}{(1+|w|)^4}\1_{[-r,r]}(|z|).
    \]
\end{lem}

Together the estimates in Lemmas \ref{knklxilemma} and \ref{offplatelemma} along with the compact support of $K_{n,k,\ell,b_I}(w,z)$ in $z$ imply
\begin{align*}
\sup_z \int \big|K_{n,k,\ell,b_I}(w,z)\big| \, dw +\sup_w \int \big|K_{n,k,\ell,b_I}(w,z)\big| \, dz\le C_\varepsilon 2^{-11k-n}.
\end{align*}
Thus
\begin{align*}
    \Big\|\sum_{I\in\mathcal{I}_J} \big(\mathrm{Id}-P_{k,A(\varepsilon),b_I,\delta_1}\big)\big[\sigma_{\ell,\varepsilon_1}\cT_{k,\ell}g_I\big]\Big\|_p
    &=\sum_{I\in\mathcal{I}_J}\sum_{n\ge 0}\Big\|\int K_{n,k,\ell,b_I}(\cdot,z) g_I(z) \, dz\Big\|_p
\end{align*}
and applying Young's inequality and the almost disjoint support of $\{g_I\}_{I\in \cI_J}$
\begin{align*}
    \sum_{I\in\mathcal{I}_J}\sum_{n\ge 0}\Big\|\int K_{n,k,\ell,b_I}(\cdot,z) g_I(z) \, dz\Big\|_p&\le  \sum_{I\in\mathcal{I}_J} \sum_{n\ge 0} C_\varepsilon 2^{-11k-n} \|g_I\|_p \\
    &\le C_\varepsilon 2^{-11k}\|g\|_p.
\end{align*}
This will complete the proof of Proposition \ref{modeldecouplingstep}.

\subsubsection{The Proof of Lemma \ref{knklxilemma}}

First, we integrate by parts in the $\xi$ variables with the differential operator 
\[
L_\xi=\langle \tfrac{w-v}{|w-v|^2}, \nabla_\xi \cdot\rangle,
\]
which will give the desired decay in $w$. Note that $\nabla_\xi \Psi=w-v$, and 
\[
|\partial_\xi^\beta\sigma_{n,2}|\le C_{|\beta|}\min\{A(\varepsilon)^{-1}\delta_1^{2}2^{k},2^{n}\}^{-|\beta|}\le C_{|\beta|} A(\varepsilon)^{|\beta|}
\]
for any multi-index $\beta$ with $|\beta|\ge 1$. Thus applying integration by parts many times with the operator $L_\xi$ gives the bound
\[
\big|(L_\xi^*)^N \sigma_{n,2}(\xi)\big|\le \frac{C_N A(\varepsilon)^N}{|w-v|^N}
\]
for any $N>0$. Since $\sigma_{n,2}$ is bounded and supported where $|\xi|\simeq 2^n$, we obtain an estimate 
\begin{equation}
\Big|\int  e^{2\pi i\Psi(w,v,z,\mu,\xi)}\sigma_{n,2}(\xi) \, d\xi\Big| \le C_N\frac{2^{3n}}{(1+A(\varepsilon)^{-1}|w-v|)^N}, \label{xiIBP}
\end{equation}
allowing us to later integrate in $w$.

By the implicit function theorem there is a constant $C_1>0$ such that if $|n-k|>C_1$ then
    \[
    |\nabla_v \Psi|=\big|-\xi+2^k\nabla_v(\mu\cdot\fS(v,z_3))\big|\ge \big||\xi|-|2^k\nabla_v(\mu\cdot \fS(v,z_3))|\big|\ge C_0\max\{2^k,2^n\}
    \]
for some $C_0>0$. We also see that $|\partial_v^\beta \Psi|\le A(\varepsilon) 2^k$ for any multi-index $\beta$ with $|\beta|\ge 2$, and $|\partial_v^\beta \sigma|\le C_{|\beta|}2^{\ell |\beta|}$ for any multi-index $\beta$ with $|\beta|\ge 1$. Since $\ell<k/2$, integrating by parts in the $v$ variables with the differential operator $L_v=\langle \tfrac{\nabla_v \Psi}{|\nabla_v \Psi|^2},\nabla_v \cdot\rangle$ gives the estimate
\[
\big|(L_v^*)^N\sigma_1(v,z,\mu)\big|\le C_N \frac{A(\varepsilon)2^\ell}{C_0\max\{2^k,2^n\}}\le C_N A(\varepsilon) \max\{2^{k/2},2^{n/2}\}^{-N}.
\]
Combining this estimate with \eqref{xiIBP}, we obtain
    \begin{align*}
    |K_{n,k,\ell,b_I}(w,z)|&\le \int\int\Big|\int e^{2\pi i \Psi(w,v,z,\mu,\xi)} \sigma_{n,2}(\xi) \, d\xi\Big| \big|(L_v^*)^{2N}\sigma_1(v,z,\mu)\big| \, dv \, d\mu  \\
    &\le C_N\int\int \frac{A(\varepsilon)}{\max\{2^k,2^n\}^{N}}\frac{1}{(1+A(\varepsilon)^{-1}|w-v|)^N} \, dv \, d\mu.
    \end{align*}
    As $\sigma_1(v,z,\mu)$ is supported where $|v|+|z|+|\mu|\le 6$ by loss of a constant depending on $\varepsilon$ we can integrate in $v$ and $\mu$ to obtain \eqref{knklxibound}.



\subsubsection{The Proof of Lemma \ref{offplatelemma}}

Suppose that $|\langle \tfrac{T_2(b_I)}{|T_2(b_I)|},\xi\rangle|\ge 3 A(\varepsilon)2^k \delta_1$. Define $\partial_{T_2(b_I)}=\langle T_2(b_I),\nabla_v\cdot\rangle$. Then by \eqref{tangentplateestimate}
\[
    |\partial_{T_2(b_I)}\Psi|\ge 2A(\varepsilon)2^k\delta_1.
\]
We can also estimate for $j\ge 1$
\[
    |\partial_{T_2(b_I)}^j \sigma_1|\le C_j A(\varepsilon)2^{\ell j},
\]
and for $j\ge 2$
\[
    |\partial_{T_2(b_I)}^j\Psi|\le C_j A(\varepsilon)2^k\le C_j A(\varepsilon)2^{\ell(j-1)}2^k\delta_1.
\]
Thus integrating by parts many times in the $T_2(b_I)$ direction and applying the estimate \eqref{xiIBP}, we obtain
\begin{align*}
    |K_{n,k,\ell,b_I}(w,z)|&\le C_{N} \int \int \frac{2^{3k}}{(1+A(\varepsilon)^{-1}|w-v|)^4} \frac{1}{(2^{k-\ell}\delta_1)^N}\, dv \, d\mu.
\end{align*}
Since $2^{k-\ell}\delta_1\ge 2^{k\varepsilon/2}$, integrating by parts in the $T_2(b_I)$ direction $\simeq 10/\varepsilon$ times and integrating over the compact support of $\sigma_1$ in $v,\mu$ gives the required estimate.
   
Next we assume that $|\langle N(b_I),\xi\rangle|\ge 3 A(\varepsilon)2^k\delta_1^2$. Define $\partial_{N(b_I)}=\langle N(b_I),\nabla_v\rangle$. Note that \eqref{normalplateestimate} implies
\begin{equation}\label{NPsilowerbound}
    |\partial_{N(b_I)}\Psi|\ge 2A(\varepsilon)2^k\delta_1^2.
\end{equation}
We claim that 
\begin{equation}\label{Nsigma1estimate}
    |\partial_{N(b_I)}^j\sigma_1(v,z,\mu)|\le C_jA(\varepsilon)\max\{2^\ell\delta_0,\varepsilon_1^{-1}\}^j
\end{equation}
for every $j\ge 1$. To see this, we use the approximation $N(b_I)=-b_Ie_1+e_3+C(b_I)b_I^2$, where $|C(b_I)|<M_0$, from \eqref{Nbestimate}. From the definition of $\sigma_1$ we see for every $j\ge 1$ and every multi-index $\beta$ with $|\beta|\le j$ that
\begin{align}
\big|(b_I\partial_{v_1})^{j-|\beta|}C(b_I)^{|\beta|}b_I^{2|\beta|}\partial_v^\beta\sigma_1(v,z,\mu)\big| &\le C_j \big(2^\ell \delta_0\big)^{j-|\beta|}\big(2^\ell\delta_0^2\big)^{|\beta|} \notag \\ 
&\le C_j \big(2^\ell\delta_0\big)^j. \label{vsigmabound}
\end{align}
Thus it suffices to check that \eqref{Nsigma1estimate} holds for mixed derivatives of the form
\[
|b_I|^{|\beta|}\partial_{v_3}^{j-|\beta|}\partial_{v'}^\beta \sigma_1(v,z,\mu),
\]
where $v'=(v_1,v_2)$, and $\beta$ is a 2-dimensional multi-index such that $|\beta|<j$. Note that
\begin{align*}
    \big||b_I|^{|\beta|}\partial_{v_3}^{j-|\beta|}\partial_{v'}^\beta \sigma_{\ell,\varepsilon_1}(v)\big|&\le C_j (2^\ell\delta_0)^{|\beta|}\varepsilon_1^{|\beta|-j} \\
    \big| |b_I|^{|\beta|}\partial_{v_3}^{j-|\beta|}\partial_{v'}^\beta \zeta(v,z,\mu)\big|&\le C_j \delta_0^{|\beta|},
\end{align*}
so it suffices to estimate 
\[
|b_I|^{|\beta|}\partial_{v_3}^{j-|\beta|}\partial_{v'}^\beta \chi_1(2^{\ell}\alpha(v,z_3)\mu\cdot\Delta(v,z_3)).
\] 
Note that terms for which no derivative hits $\mu\cdot\Delta(v,z_3)$ will be negligible since $|\mu\cdot\Delta(v,z_3)|\simeq 2^{-\ell}$.  
Using \eqref{delta1approx}, \eqref{blowdowndelta2}, and a Taylor expansion about $z_3=0$ we see that
\begin{align*}
    \big||b_I|^{|\beta|}\partial_{v_3}^{j-|\beta|}\partial_{v'}^\beta \big[\mu\cdot \Delta(v,z_3)\big]\big|&\le \delta_0^{|\beta|}\partial_{v'}^\beta (\mu\cdot\Delta_{v_3^j}(v,0))+A(\varepsilon)\delta_0^{|\beta|+1} \\
    &=A(\varepsilon)\delta_0^{|\beta|+1}.
\end{align*}
Thus we see  by differentiation of compositions and products
\begin{equation}
   \big| |b_I|^{|\beta|}\partial_{v_3}^{j-|\beta|}\partial_{v'}^\beta \sigma_1(v,z,\mu)\big|\le C_j C_jA(\varepsilon) \max\{2^\ell\delta_0,\varepsilon_1^{-1}\}^j \label{v3sigmabound},
\end{equation}
and by combining \eqref{vsigmabound} and \eqref{v3sigmabound} the claim \eqref{Nsigma1estimate} is proven.

To integrate by parts we also need to show that for $j\ge 2$
\[
    |\partial_{N(b_I)}^j\Psi|\le C_jA(\varepsilon)2^k\max\{2^{\ell}\delta_0,\varepsilon_1^{-1}\}^{j-1}\delta_1^{2}.
\]
In fact, we claim that
\begin{equation}
    |\partial_{N(b_I)}^j\Psi|\le C_jA(\varepsilon)2^k\delta_0^2 \label{NPsibound}
\end{equation}
for $j\ge 2$.
We use \eqref{Nbestimate} again to see that
\[
    \partial_{N(b_I)}^j\Psi=\langle b_Ie_1+e_3+C(b_I)b_I^2, \nabla_v\rangle^j\Psi
\]
where again $|C(b_I)|\le M_0$. Rearranging terms using the fact that $|\partial_{v}^\beta\Psi|\le A(\varepsilon)2^k$ for any multi-index $\beta$ with $|\beta|\ge 2$ we obtain
\[
\partial_{N(b_I)}^j\Psi=b_I\Psi_{v_1v_3^{j-1}}+\Psi_{v_3^j}+A(\varepsilon)2^k b_I^2R_{10}(v,z_3).
\]
Next, we estimate via a Taylor expansion about $z_3=0$,
\begin{align*}
    \Psi_{v_1v_3^{j-1}}(v,z_3)&=\mu\cdot\fS_{v_1v_3^{j-1}}(v,0)+2M_jz_3R_{11}(v,z_3) \\
    \Psi_{v_3^j}(v,z_3)&=\mu\cdot\fS_{v_3^j}(v,0)+z_3\mu\cdot\fS_{v_3^jz_3}(v,0)+2M_jz_3^2 R_{12}(v,z_3).
\end{align*}
From \eqref{S1wderivs}, \eqref{S2wderivs}, and \eqref{S1zwderivs} we see that
\begin{align*}
    \mu\cdot \fS_{v_1v_3^{j-1}}(v,0)&=0, \\
    \mu\cdot \fS_{v_3^j}(v,0)&=0, \\
    \fS^1_{v_3^j z_3}(v,0)&=0.
\end{align*}
Moreover \eqref{blowdownS2} ensures that
\[
    \fS^2_{v_3^j z_3}(v,0)=0, \ j\ge 2.
\]
Hence for $j\ge 2$
\[
    |\partial_{N(b_I)}^j\Psi|\le C_j A(\varepsilon)\delta_0^2\le C_j A(\varepsilon)\varepsilon_1^{-j+1}\delta_1^2,
\]
satisfying the claim \eqref{NPsibound}. 

Now that we have verified the conditions \eqref{NPsilowerbound}, \eqref{Nsigma1estimate}, and \eqref{NPsibound}, we integrate by parts $M$ times in the $N(b_I)$ direction and apply \eqref{xiIBP} to obtain
\begin{align*}
    |K_{n,k,\ell,b_I}(w,z)|&\le\Big(\frac{1}{\min\{2^{k-\ell}\delta_1^2/\delta_0,2^k \delta_1^2\varepsilon_1\}}\Big)^{M} \int\int \frac{1}{(1+A(\varepsilon)^{-1}|w-v|)4} dv \, d\mu.
\end{align*}
Since $\delta_1\ge 2^{-\ell(1-\varepsilon/2)}$ and $\delta_1\ge 2^{-\ell\varepsilon/4}\delta_0$, we have
\begin{align*}
    2^{k-\ell}\delta_1^2/\delta_0&\ge 2^{k-2\ell+\ell\varepsilon/4}\ge 2^{k\varepsilon/8} \\
    2^k\delta_1^2\varepsilon_1&\ge  2^{k-2\ell+\ell\varepsilon/2}\ge 2^{k\varepsilon/4}.
\end{align*}
So if $M\simeq 50/\varepsilon$ and we integrate over the compact support of $\sigma_1$ in $v$ and $\mu$ we obtain the desired estimate.

\subsection{Families of changes of variables}
We use the family of changes of variables used in \cite{PrSe19} to reduce the general case to the model case. Let $P^\circ=(a^\circ,y^\circ)\in\cM$, with $y^{\circ}=(S^1(a^\circ,b^\circ),S^2(a^\circ,b^\circ),b^\circ)$. For $r>0,q>0$ let 
\[
Q(r)=\{(x_1,x_2,x_3) \ : \ |x-a^\circ|\le  r\}
\]
and
\[
I(r)=\{y_3 \ : \ |y_3-b^\circ|\le r\}.
\]
For $i=1,2$, let $S^i$ be smooth functions in a neighborhood of $Q(2r_0)\times I(2r_0)$, for some $r_0>0$. We assume that $\Delta_1(x,y_3)=\det(S^1_x, S^2_x, S^1_{xy_3})\ne 0$ on $Q(2r_0)\times I(2r_0)$. Choose $M>0$ so that
\[
M>2+\|S\|_{C^5(Q(2r_0)\times I(2r_0))}+\max_{(x,y_3)\in Q(2r_0)\times I(2r_0)} |\Delta_1(x,y_3)|^{-1}
\]
For $a\in Q(2r_0)$, $b\in I(2r_0)$, let 
\begin{align*}
\Gamma_1(x,y_3)&=\det(S_x^1, \, S_{x y_3}^2, \, S_{xy_3}^1) \\
\Gamma_2(x,y_3)&=\det(S_{x y_3}^1, \, S_{x}^2, \, S_{xy_3}^2),
\end{align*}
and let $\rho(a,b)\in\rr^3$ be defined by
\[
(\rho_1,\rho_2,\rho_3):=\frac{1}{\Delta_1(a,b)}(-\Gamma_2(a,b),\Gamma_1(a,b),\Delta_2(a,b)).
\]
For $(x,y_3),(a,y_3)\in Q(r_0)\times I(2r_0)$, consider the map
\[
(x,a,y_3)\mapsto \fw(x,a,y_3)\in\rr^3
\]
given by
\begin{align*}
    \fw_1(x,a,b)&=S^1(x,b)-S^1(a,b) \\
    \fw_2(x,a,b)&=S^2(x,b)-\rho_3(a,b)S^1(x,b)-S^2(a,b)+\rho_3(a,b)S^1(a,b) \\
    \fw_3(x,a,b)&=S_{y_3}^1(x,b)-S_{y_3}^1(a,b).
\end{align*}
Then
\[
\det(D\fw/Dx)=\det(S_{x}^1(x,b),\, S_x^2-\rho_3(a,b)S_x^1(x,b),\, S_{xy_3}^1(x,b))=\Delta^1(x,b)\ne 0.
\]
By the implicit function theorem, there exists $r_1>0$ with $r_1<r_0$ such that $|w|_\infty<2r_1$, $a\in Q(2r_1)$, $b\in I(2r_1)$, the equation $\fw(x,a,b)=w$ is solved by a unique $C^\infty$ function $x=\fx(w,a,b)$ such that for $|x|_\infty<(50M^5)^{-1}r_1$,  $\fx(\fw(x,a,b),a,b)=x.$

We also change variables in $y$. Define $\fz:\rr^2\times Q(2r_0)\times I(2r_0)^2\to \rr^3$ by
\begin{align*}
    \fz_1(y,a,b)&=y_1-S^1(a,y_3) \\
    \fz_2(y,a,b)&=y_2-\rho_3(a,b)y_1-S^2(a,y_3)+\rho_3(a,b)S^1(a,y_3)-(y_3-b)\sum_{i=1}^2 \rho_i(y_i-S^i(a,y_3)) \\
    \fz_3(y,a,b)&=y_3-b.
\end{align*}
For $z_3,b\in I(r_3)$, $a\in Q(2r_0)$, where $r_3<\min\{r_1,(24M^4)^{-1}\}$, we define the inverse $z\mapsto \fy(z,a,b)$ by
\begin{align*}
    \fy_1(z,a,b)&=z_1+S^1(a,b+z_3) \\
    \fy_2(z,a,b)&=\frac{z_2+z_1(\rho_3(a,b)+\rho_1(a,b)z_3)+(1-z_3)S^2(a,b+z_3)}{1-\rho_2(a,b)z_3} \\
    \fy_3(z,a,b)&=b+z_3.
\end{align*}


\begin{lem}[\cite{PrSe19}]
    The function $\fx,\fy$ defined above have the following properties. 
    \begin{enumerate}
        \item $\fx(0,a,b)=a$, $\fy(0,a,b)=(S^1(a,b),S^2(a,b),b)$, $\fy_3(z,a,b)=b+z_3$. 
        \item $\det\Big(\tfrac{D\fx(w,a,b)}{dw}\Big)=\tfrac{1}{\Delta^1(\fx(w,a,b),y_3)}$.
        \item Let $\rho=\rho(a,b)$, and let \[
        B(z_3,a,b)=\Big(\begin{smallmatrix} 1 & 0 \\ -\rho_3-\rho_1 z_3 & 1-\rho_2 z_3\end{smallmatrix}\Big).
        \]
        Then for $|z_3|\le r_3$, $a\in Q(2r_1,2q_1)$, $|w|\le r_1$
        \[
        B(z_3,a,b)\Big(\begin{smallmatrix} S^1(\fx(w,a,b),b+z_3)-\fy_1(z,a,b) \\ S^2(\fx(w,a,b),b+z_3)-\fy_2(z,a,b) \end{smallmatrix} \Big)=\Big(\begin{smallmatrix} \fS^1(w,z_3,a,b)-z_1 \\ \fS^2(w,z_3,a,b)-z_2\end{smallmatrix}\Big),
        \]
        where $\fS^i$ are $C^\infty$ with
        \[
        (\fS^1,\fS^2,\fS^1_{z_3})\big|_{(w,0,a,b)}=w
        \] 
        and $\fS^2_{wz_3}(0,0,a,b)=0.$
        \item Let
        \begin{align*}
        \Delta^i_S(x,y_3)&=\det(S^1_x,S^2_x,S^i_{xy_3})\big|_{(x,y_3)} \\
        \Delta^i_\fS(x,y_3,a,b)&=\det(\fS^1_w, \fS^2_w, \fS^i_{wz_3})\big|_{(w,z_3,a,b)}.
        \end{align*}
        Then for $(\tau_1,\tau_2)=(\mu_1,\mu_2)B(z_3,a,b),$
        \[
        \sum_{i=1}^2 \tau_i \Delta^i_S(\fx(w,a,b),b+z_3)=\frac{\Delta^1_S(\fx(w,a,b),b)}{1-\rho_2(a,b)z_3}\sum_{i=1}^2 \mu_i\Delta^i_\fS(w,z_3,a,b).
        \]
        \item Let \[
        \kappa(a,b)=\Gamma_2\Delta^1_S-\Gamma_1\Delta^2_S+\Delta^1_S\Delta^2_{S,y_3}-\Delta^2_S\Delta^1_{S,y_3}\Big|_{(a,b)}.
        \]
        Then
        \[
        \fS^2_{w_3z_3z_3}(0,0,a,b)=\frac{\kappa(a,b)}{\Delta^1_S(a,b)}.
        \]
    \end{enumerate}
\end{lem}

The proof is found in \cite{PrSe19}.

\subsection{Decoupling in the General Case}

\begin{prop}\label{generaldecouplingstep}
Let $0<\varepsilon<\tfrac1{10}$, $k\gg 1$, $\ell\le k/2$. Let $\delta_0\in (2^{-\ell(1-\varepsilon)},2^{-\ell\varepsilon})$, and define $0<\delta_1<\delta_0$ such that 
\[
\max\{2^{-\ell(1-\varepsilon/2)},\delta_0 2^{-\ell\varepsilon/4}\}<\delta_1<\delta_0.
\]
Define $\varepsilon_1=(\delta_1/\delta_0)^2$. Let $J$ be an interval of length $\delta_0$ within $r_0$ of $b^\circ$, and let $\cI_J$ be a collection of intervals which have disjoint interior, intersecting $J$. For each $I\in\cI_J$, define $f_I(y)=f(y)\1_{I}(y_3)$. Then for $2\le p\le 6$ 
\[
\norm{\sum_{I\in\cI_J} \cR_{k,\ell} f_I}_p\le C_\varepsilon (\delta_0/\delta_1)^{1/2-1/p+\varepsilon}\Big(\sum_{I\in\cI_J} \norm{\cR_{k,\ell}f_I}_p^p\Big)^{1/p}+C_{N,\varepsilon}2^{-kN}\|f\|_p.
\]
\end{prop}

\bp
Fix $b\in J$. Let $\sigma_0\in C_c^\infty$ supported in $(-1,1)$ such that $\sigma_0\ge 0$ everywhere and $\sum_{n\in\zz}\sigma_0(\cdot-n)=1$. For $n\in\zz$ define $\sigma_{n}(x)=\sigma_0(x-n)$ and for $a\in\zz^3$ let 
\[
\varsigma_{\ell,\varepsilon_1}(x,a,b)=\sigma_{a_1}(2^{\ell}\fw_1(x,a,b))\sigma_{a_2}(2^\ell\fw_2(x,a,b))\sigma_{a_3}(\varepsilon_1^{-1}\fw_3(x,a,b)).
\]
Then since $|x|,|a|<r_2/2$ and $|b|<r_3/2$, we have that
\[
\varsigma_{\ell,\varepsilon_1}(\fx(w,a,b),a,b)=\sigma_{a_1}(2^\ell w_1)\sigma_{a_2}(2^\ell w_2)\sigma_{a_3}(\varepsilon_1^{-1}w_3)=:\sigma_{\ell,\varepsilon_1}(w,a),
\]
and $\sum_{a\in\zz^3} \sigma_{\ell,\varepsilon_1}(w,a)=1$ with finite overlap for all $|w|<r_2/2$. Then by H\"older's inequality
\[
\norm{\sum_{I\in\cI_J}\cR_{k,\ell}f_I}_p\le C_p \Big(\sum_{a\in\zz^3}  \norm{\sum_{I\in\cI_J} \varsigma_{\ell,\varepsilon_1,a,b}\cR_{k,\ell} f_I}_p^p\Big)^{1/p}.
\]
Note that the terms vanish for $|a|>r_2$. Fix $a$. Write $g(z,a,b)=f(\fy(z,a,b))$. Apply changes of variables $y=\fy(z,a,b)$ and $\tau={B^\intercal}^{-1}(z_3,a,b)\mu$, noting that 
\[
\det(D\fy/dz)\det B=1
\]
so that
\[
\cR_{k,\ell}f(x)=2^{2k}\iint e^{i2^k \langle \mu,\fS(\fw(x,a,b)-z'\rangle} \tilde{\chi}_{k,\ell}(x,z,\mu,a,b)g(z,a,b) \, d\mu dz,
\]
with 
\begin{align*}
\tilde{\chi}_{k,\ell}(x,z,\mu,a,b)&= \chi(x,\fy(z,a,b))\eta_1(|{B^\intercal}^{-1}(z_3,a,b)\mu|)\eta\Big(2^\ell\tfrac{\Delta^1(x)}{1-\rho_3(a,b)z_3} \\
&\qquad \times(\mu_1\Delta^1_{\fS}(\fw(x,a,b),z_3,a,b)+\mu_2\Delta^2_{\fS}(\fw(x,a,b),z_3,a,b))\Big).
\end{align*}
Thus we see that
\[
\varsigma_{\ell,\varepsilon_1,a,b}(\fx(w,a,b))\sum_{I\in\cI_J} \cR_{k,\ell} f_I(\fx(w,a,b))= \sigma_{\ell,\varepsilon_1,a,b}(w)\sum_{I\in\cI_J} T_{k,\ell,a,b}g_I(w),
\]
where $g_I(z,a,b)=g(z,a,b)\1_{-b+I}(z_3)$ and $\cT_{k,\ell,a,b}\equiv \cT_{k,\ell}$ from the model case. Define
 \begin{align*}
     M_n(a,b)&\ge 2+||\fS^1(\cdot,a,b)||_{C^{n+5}([-r_0,r_0]^4)}+||\fS^2(\cdot,a,b)||_{C^{n+5}([-r_0,r_0]^4)} \\
    \tilde{A}(\varepsilon)&=\sup_{a,b\in[-r_0,r_0]^4} \max\{3^{\lceil3/\varepsilon\rceil+2}M_{\lceil 3/\varepsilon\rceil}(a,b),\kappa_0(a,b)(1+4^{\lceil3/\varepsilon\rceil+2}M_{\lceil 3/\varepsilon\rceil}^4(a,b))\};
 \end{align*}
these are the uniform versions of \eqref{Mn} and \eqref{Aepsilon} respectively. We can then write
\begin{align*}
\norm{\varsigma_{\ell,\varepsilon_1,a,b}\sum_{I\in\cI_J}\cR_{k,\ell}f_I}_p&=\Big(\int \Big|\varsigma_{\ell,\varepsilon_1,a,b}(\fx(w,a,b))\sum_{I\in\cI_J}\cR_{k,\ell}f_I(\fx(w,a,b))\|^p |\det(\tfrac{D\fx}{Dw})| \, dw \Big)^{1/p} \\
&\lesssim \norm{\sigma_{\ell,\varepsilon_1,a,b}\sum_{I\in\cI_J} \cT_{k,\ell} g_I}_p
\end{align*}
by the uniform upper bound on $|\det(\tfrac{D\fx}{Dw})|$. Then we can apply Proposition \ref{modeldecouplingstep} with $A(\varepsilon)=\tilde{A}(\varepsilon)$ to get
\begin{align*}
\norm{\sigma_{\ell,\varepsilon_1,a,b}\sum_{I\in\cI_J} \cT_{k,\ell} g_I}_p&\le C_\varepsilon(\delta_0/\delta_1)^{1/2-1/p+\varepsilon}\Big(\sum_{I\in\cI_J}\big\|\sigma_{\ell,\varepsilon_1,a,b} \cT_{k,\ell}g_I\big\|_p^p\Big)^{1/p} \\
& \qquad +C_\varepsilon 2^{-10k}2^{-2\ell} \varepsilon_1\|g\|_p.
\end{align*}
Then undoing the changes of variables above (and applying the uniform lower bounds on $|\det(\tfrac{D\fx}{Dw})|$) we may bound this by
\begin{align*}
    C'_\varepsilon(\delta_0/\delta_1)^{1/2-1/p+\varepsilon}\Big(\sum_{I\in\cI_J}\big\|\varsigma_{\ell,\varepsilon_1,a,b} \cR_{k,\ell}f_I\big\|_p^p\Big)^{1/p}+C_\varepsilon 2^{-10k}2^{-2\ell} \varepsilon_1\|f\|_p.
\end{align*}
Finally, we recombine our partition of unity in $x$ using the fact that there are at most $C2^{2\ell}\varepsilon_1^{-1}$ many $a\in\zz^3$ for which $\sigma_{\ell,\varepsilon_1,a,b}$ is nonzero, to get
\begin{align*}
\norm{\sum_{I\in\cI_J}\cR_{k,\ell}f_I}_p&\le C_p \Big(\sum_{\substack{a\in\zz^3 \\ |a|_\infty<r_2}}  \norm{\sum_{I\in\cI_J} \varsigma_{\ell,\varepsilon_1,a,b}\cR_{k,\ell} f_I}_p^p\Big)^{1/p} \\
&\le C_p C_\varepsilon(\delta_0/\delta_1)^{1/2-1/p+\varepsilon}\Big(\sum_{\substack{a\in\zz^3 \\ |a|_\infty<r_2}}\sum_{I\in\cI_J}\big\|\varsigma_{\ell,\varepsilon_1,a,b} \cR_{k,\ell}f_I\big\|_p^p\Big)^{1/p} \\
&\qquad +\sum_{\substack{a\in\zz^3 \\ |a|_\infty<r_2}} C_\varepsilon 2^{-2\ell}\varepsilon_1 2^{-10k}\|f\|_p \\
&\le C_\varepsilon(\delta_0/\delta_1)^{1/2-1/p+\varepsilon}\Big(\sum_{I\in\cI_J} \|\cR_{k,\ell}f_I\|_p^p\Big)^{1/p} +C_\varepsilon 2^{-10k}\|f\|_p.
\end{align*}

\ep

\subsection{Iteration of the Decoupling Step}

Let $\delta_0=2^{-\ell\varepsilon}$, and define $\delta_j=\delta_{j-1}2^{-\ell\varepsilon/4}$ for $j=1,2,...$ Note that this implies $\varepsilon_1=(\delta_1/\delta_0)^2=2^{-\ell\varepsilon/2}$. We will iterate the estimate in Proposition \ref{generaldecouplingstep} until $\delta_j\le 2^{-\ell(1-\varepsilon)}$. Let $j^*$ be the smallest $j$ such that $\delta_j<2^{-\ell(1-\varepsilon)}$. Clearly $j^*\lesssim 1/\varepsilon$ and $2^{-\ell(1-\varepsilon/2)}\le \delta_{j^*}\le 2^{-\ell(1-\varepsilon)}$. 

For $j=0,1,2,...$ let $I_j$ denote an interval of length $\delta_j$ inside $[b^\circ-r_0,b^\circ+r_0]$, and let $\cI_{I_j}$ denote the collection of intervals $I_{j+1}$ of length $\delta_{j+1}$ intersecting $I_j$ with disjoint interior. Finally, let $J=[b^\circ-r_0/2,b^\circ+r_0/2]$ and let $\cI_{J,j}$ denote the collection of intervals $I_{j}$ of length $\delta_{j}$ intersecting $J$ with disjoint interiors. Then since $r_0<1$ and $\delta_0=2^{-\ell\varepsilon}$, using H\"older's and Minkowski's inequalities we have
\begin{equation}\label{basecase}
\|\cR_{k,\ell} f\|_p\lesssim 2^{\ell\varepsilon/p'}\Big(\sum_{I_0\in \cI_{J,0}} \norm{\cR_{k,\ell} f_{I_0}}_p^p\Big)^{1/p}.
\end{equation}
The function and operator $\cR_{k,\ell} f_{I_0}$ now satisfy the conditions of Proposition \ref{generaldecouplingstep}. We claim that for each $0\le j\le j^*$,
\begin{align} \label{induction}
\|\cR_{k,\ell} f\|_{p}&\lesssim C(\varepsilon)^j  2^{\ell\varepsilon/(p')} (\delta_0/\delta_j)^{1/2-1/p+\varepsilon}\Big(\sum_{I_j\in \cI_{J,j}} \|\cR_{k,\ell} f_{I_j}\|_p^p\Big)^{1/p} \\
& \qquad  + j 2^{2\ell} C(\varepsilon)^{j} 2^{-10k}\|f\|_p.\notag
\end{align}
The case $j=0$ follows immediately from \eqref{basecase}. Assume \eqref{induction} holds for some $j$. Then by applying Proposition \ref{generaldecouplingstep} we get
\begin{align}
\Big(\sum_{I_j\in\cI_{J,j}}\| \cR_{k,\ell} f_{I_j}\|_{p}^p\Big)^{1/p}&\le \Big(\sum_{I_j\in\cI_{J,j}} \Big[C(\varepsilon)\big(\tfrac{\delta_j}{\delta_{j+1}}\big)^{1/2-1/p+\varepsilon}\Big(\sum_{I_{j+1}\in I_{I_j}} \|\cR_{k,\ell} f_{I_{j+1}}\|_p^p\Big)^{1/p} \notag \\
& \qquad + C(\varepsilon) 2^{-10k} \|f\|_p^p\Big)^{1/p}\Big]^p\Big)^{1/p}\label{applydecoupling} \\
&\le  C(\varepsilon)\big(\tfrac{\delta_j}{\delta_{j+1}}\big)^{1/2-1/p+\varepsilon}\Big(\sum_{I_{j+1}\in\cI_{J,j+1}}\|\cR_{k,\ell} f_{I_{j+1}}\|_p^p\Big)^{1/p} \notag \\
& \qquad + C(\varepsilon) \delta_{j}^{-1/p}2^{-10k}\|f\|_p. \notag
\end{align}
Plugging the above estimate into \eqref{induction} gives us
\begin{align*}
\|\cR_{k,\ell} f\|_p &\le C(\varepsilon)^{j+1}2^{\ell\varepsilon/(p')}\big(\tfrac{\delta_{0}}{\delta_{j+1}}\big)^{1/2-1/p+\varepsilon}\Big(\sum_{I_{j+1}\in\cI_{J,j+1}} \|\cR_{k,\ell} f_{I_{j+1}}\|_p^p\Big)^{1/p} \\ 
& \qquad + C(\varepsilon)^{j}2^{\ell\varepsilon/(p')}\big(\tfrac{\delta_0}{\delta_j}\big)^{1/2-1/p+\varepsilon}C(\varepsilon)\delta_j^{-1/p}2^{-10k}\|f\|_p \\
& \qquad + j2^{2\ell}C(\varepsilon)^{j}2^{-10k}\|f\|_p.
\end{align*}
Using the fact that $\delta_0=2^{-\ell\varepsilon}$, $\delta_{j}\ge 2^{\ell(1-\varepsilon/2)}$ for $j\le j^*$, and $2\le p\le 6$, the last two terms of the above inequality are bounded by
\[
(j+1)C(\varepsilon)^{j+1}2^{2\ell} 2^{-10k}\|f\|_p,
\]
proving the claim. 

We apply \eqref{induction} for $j=j^*$ and use the fact that $j^*\le 4/\varepsilon$ as well as the assertion
\[
\frac{\varepsilon}{p'}-\frac{\varepsilon}{2}+\frac{\varepsilon}{p}-\varepsilon^2-\frac{\varepsilon}{4}+\frac{\varepsilon}{2p}+\frac{\varepsilon^2}{2}\le \varepsilon
\]
to deduce
\begin{align}
\|\cR_{k,\ell} f\|_p & \le  C(\varepsilon)^{4/\varepsilon}2^{\ell\varepsilon/p'}2^{-\ell\varepsilon\big(\frac12-\frac1p+\varepsilon\big)}2^{\ell\big(1-\frac{\varepsilon}{2}\big)\big(\frac12-\frac1p+\varepsilon\big)}\Big(\sum_{I_{j^*}\in \cI_{J,j^*}}\|\cR_{k,\ell} f_{I_{j^*}}\|_p^p\Big)^{1/p} \\
& \qquad + \tfrac{4}{\varepsilon}C(\varepsilon)^{4/\varepsilon} 2^{-10k+2\ell}\|f\|_p  \notag \\
& \lesssim_{\varepsilon} 2^{\ell(1/2-1/p+2\varepsilon)}\Big(\sum_{I_{j^*}\in \cI_{J,j^*}} \|\cR_{k,\ell}f_{I_{j^*}}\|_p^p\Big)^{1/p}+ C(\varepsilon)2^{-9k}\|f\|_p.
\end{align}
Picking $\varepsilon'=2\varepsilon$ completes the proof.

\section{$L^2$ Estimates}\label{L2section}
The methods in this section draw from the work of Comech in \cite{Co97,Co99}, which was itself influenced by the estimates proven in \cite{PhSt91}. While Comech proved $L^2$ regularity estimates for fold and finite type conditions, here we prove $L^2$ estimates for a general class of oscillatory integral operators associated to fold blowdown singularities in $d$ dimensions. Let
\[
\cA_k f(x)=\int e^{i2^k\phi(x,y)}f(y)\sigma(x,y) dy,
\]
where $x,y\in\rr^d$, $\phi\in C^\infty(\rr^d\times\rr^d)$, and $\sigma\in C_c^\infty(\rr^d\times\rr^d)$. The canonical relation associated to this oscillatory integral operator is given by
\[
\{(x,\phi_x)\times (y,\phi_y) \ : \ x\in\rr^d, \ y\in\rr^d\}
\]
We write the projections $\pi_L:(x,y)\mapsto (x,\phi_x)$ and $\pi_R:(x,y)\mapsto (y,\phi_y)$. The projections are degenerate on the variety $\cL$ where the determinant of the mixed Hessian of $\phi$ vanishes. Let $h(x,y)=\det \phi_{xy}$. We assume that $\pi_L$ is a fold and $\pi_R$ is a blowdown on $\cL$. We may choose the support of $\sigma$ small enough and choose coordinates $x=(x',x_d)$, $y=(y',y_d)$ in $\rr^{d-1}\times\rr$ vanishing at a reference point $P^\circ=(x^\circ,y^\circ)$ so that
\[
\phi_{x'y'}(P^\circ)=I_{d-1}, \quad \phi_{x_dy'}(P^\circ)=0, \quad \phi_{x'y_d}(P^\circ)=0, \quad \phi_{x_dy_d}(P^\circ)=0.
\]
Let $\phi^{x'y'}=\phi_{x'y'}^{-1}$, and we can define the kernel fields
\begin{align*}
  V_R&=\partial_{x_d}-\phi_{x_dy'}(\phi^{x'y'})^{\intercal}\partial_{x'} \\
  V_L&=\partial_{y_d}-\phi_{x'y_d}\phi^{x'y'}\partial_{y'}.
\end{align*}
By the assumption on $\pi_L$, in other words that $\phi_{xy}$ has corank at most 1 and
\[
h(x,y)=0 \implies |V_L h(x,y)|\ge c_L>0.
\]
Since we assume that $V_R$ is a blowdown, i.e. that $V_R$ is tangent to the singularity surface $\cL$, we see that
\[
h(x,y)=0 \implies V_R^j h(x,y)=0 \ \forall j\ge 0.
\]
Assuming small enough support of $\sigma$ we may assume that for $(x,y)\in\supp( \sigma)$
\[
\max\{|\phi_{x'y_d}(x,y)|,|\phi_{x_d y'}(x,y)|\}<\varepsilon.
\]
Note that this implies
\[
|(V_L-\partial_{y_d}) h(x,y)|\le \varepsilon\|\phi\|_{C^3}.
\]
We decompose first in distance to the singularity surface $\cL$; for $\ell\le \ell_0=\lfloor \tfrac{k}{2+\varepsilon}\rfloor<\tfrac{k}{2}$, we define
\[
\cA_{k,\ell}f(x):=\int e^{i2^k\phi(x,y)}f(y)\sigma(x,y)\chi(2^\ell h(x,y)) dy.
\]
\begin{thm}\label{L2general}
    For $\varepsilon>0$ and $\ell\le \ell_0$
    \[
    \|\cA_{k,\ell}f\|_2\le C_\varepsilon 2^{(\ell-dk)/2+\ell\varepsilon}\|f\|_2.
    \]
\end{thm}

To prove Proposition \ref{L2boundA} for $\cR_{k,\ell}$ we apply a partial Fourier transform (in the $y'$ variables) then change variables $2^k\mu=\tau$, which satisfies the conditions for Theorem \ref{L2general} with $d=3$. 

\subsection{Proof of Theorem \ref{L2general}}
First we note that if Theorem \ref{L2general} holds for $\ell<\ell_0$ then the global estimate 
\[
\Big\|\sum_{\ell\le \ell_0}\cA_{k,\ell}f\Big\|_2\lesssim 2^{\tfrac{k}{4}-\tfrac{dk}{2}}\|f\|_2.
\]
from \cite{GrSe94} implies the result for $\ell=\ell_0$. 

Since $k=(2+\varepsilon)\ell_0$, we see $\tfrac{k}{4}\le \tfrac{\ell}{2}+\varepsilon\tfrac{\ell}{2}$, hence by triangle inequality
\[
\|\cA_{k,\ell_0} f\|_2\lesssim 2^{(\ell-dk)/2+\ell\varepsilon}\|f\|_2.
\] 

For $\ell<\ell_0$ we decompose our operator further and use methods of the proof of the Calderon-Vaillancourt theorem, following the ideas of Comech in \cite{Co97}. We decompose our operator along small boxes in $y$-space, by way of cutoffs
\[
\chi_{\vec{m}}(y)=\prod_{j=1}^d\chi(2^\ell y_j-m_j)).
\]




We fix $k,\ell$ for now and let $\cA_{m}:=\cA_{k,\ell}[\chi_m(y)\cdot]$. Then $\cA_{m}\cA^*_{\tilde{m}}$ has Schwartz kernel 
\begin{equation*}
    K^{\cA\cA^*}_{m,\tilde{m}}(x,w)= \int e^{i2^k(\phi(x,y)-\phi(w,y))} \sigma(x,w,y,k,\ell)\enskip dy,
\end{equation*}
where the amplitude is given by
\begin{align*}
    \sigma(x,w,y)&= \chi(2^\ell h(x,y)) |\chi_m(y)|^2\overline{\chi(2^{\ell} h(w,y))}. 
\end{align*}

Similarly, the Schwartz kernel for $\cA^*_{m}\cA_{\tilde{m}}$ is given by
\begin{equation*}
    K^{\cA^*\cA}_{m,\tilde{m}}(y,z)=\int e^{i2^k(\phi(x,y)-\phi(x,z))} \tilde{\sigma}(x,y,z) \enskip dx,
\end{equation*}
where 
\begin{align*}
    \tilde{\sigma}(x,y,z)&= \chi(2^{\ell} h(x,y))\chi_m(y)\overline{\chi(2^\ell h(x,z)}\overline{\chi_{\tilde{m}}(z)}.
\end{align*}
By splitting our operator $\cA_{k,\ell}$ into a finite number of collections of $\{\cA_{m}\}$ we may assume that if $m_j\ne \tilde{m}_j$ then $|m_j-\tilde{m}_j|> \max\{\tfrac{15}{c_L},2\sqrt{d}\}$. We may also assume through the loss of a constant $2^{\ell\varepsilon}$ that the kernels of $\cA_m$ are supported 
where $|x_d|\le 2^{-\ell\varepsilon}$.

We first prove the following lemmas.

\begin{lem}\label{individualestimate}
     There exists a constant $C>0$ such that 
     \[
         \|\cA_{m}\|_{L^2\to L^2}\le C 2^{\ell-dk}.
     \]
 \end{lem}

\begin{lem}\label{almostorth} For any $N>0$, and $\ell<\ell_0=\lfloor \tfrac{k}{2+\varepsilon}\rfloor$ the following estimates hold.
    \be
        \item[(a)] If $m\ne\tilde{m}$ then
        \[
            \|\cA_{m}\cA^*_{\tilde{m}}\|=0.
        \]
        \item[(b)] If $m\ne \tilde{m}$ and $|m'-\tilde{m}'|\le \tfrac{c_L}{10\|\phi\|_{C^3}} |m_d-\tilde{m}_d|$ then
        \[
            \|\cA^*_{m}\cA_{\tilde{m}}\|=0.
        \]
        \item[(c)] If $m\ne \tilde{m}$ and $|m'-\tilde{m}'|\ge \tfrac{c_L}{10\|\phi\|_{C^3}} |m_d-\tilde{m}_d|$ then
        \[
            \|\cA^*_{m}\cA_{\tilde{m}}\|\lesssim_N 2^{\ell-dk}\big(2^{k-2\ell}|m-\tilde{m}|\big)^{-N}.
        \]
    \ee
\end{lem}
A few remarks. First, the estimates in Lemma \ref{almostorth} do not rely on the blowdown assumption and essentially reproves the results of Comech in \cite{Co97}, albeit through a slightly different approach. Second, the separation of $\ell$ from $k/2$ is necessary for the proof of Lemma \ref{almostorth}, but not Lemma \ref{individualestimate}.

\subsubsection{Proof of Lemma \ref{individualestimate}} 
Since $|\phi_{x'y'}|>c>0$ the set of equations $\nabla_{y'}(\phi(x,y)-\phi(w,y))=0$ is solved uniquely by $x'=\fx'(w,x_d,y)$. By the implicit function theorem we can see that 
\[
\frac14|x'-\fx'(w,x_d,y)|\le |\phi_{y'}(x,y)-\phi_{y'}(w,y)|\le 4|x'-\fx'(w,x_d,y)|.
\]
A further set of calculations reveal that
\begin{align*}
\phi_{y_d}(\fx'(w,x_d,y),x_d,y)-\phi_{y_d}(w,y)&=\sum_{j=0}^N V_R^j [\det \phi_{xy} (\det \phi_{x'y'})^{-1}](w,y) \frac{(x_d-w_d)^{j+1}}{(j+1)!} \\
&\qquad +b(w,y)(x_d-w_d)^{N+2},
\end{align*}
where $b$ is uniformly bounded and $N$ is chosen such that $|x_d-w_d|^N\le 2^{-\ell}$. Since $\pi_L$ is a fold and $V_L\big|_{(0,0)}=\partial_{y_d}$, we see that $h(x,y)=0$ is solved uniquely by $y_d=\fy_d(x,y')$ near $0$. From this,
\[
\frac14 |y_d-\fy_d(x,y')|\le |h(x,y)|\le 4|y_d-\fy_d(x,y')|.
\]
Because $\pi_R$ is a blowdown and the bounds on $h$
we see that
\[
|V_R^jh(x,y)|=|V_R^jh(x,y',\fy_d(x,y'))+(y_d-\fy_d(x,y'))\partial_{y_d}V_R^jh(x,y',z_d)|\le C 2^{-\ell}
\]
implying by the properties of differentiation of products
\[
|\phi_{y_d}(\fx'(w,x_d,y),x_d,y)-\phi_{y_d}(w,y)|\ge c 2^{-\ell} |x_d-w_d|.
\]
Thus 
\begin{align*}
|\phi_{y_d}(x,y)-\phi_{y_d}(w,y)|&\ge |\phi_{y_d}(x,y)-\phi_{y_d}(\fx'(w,x_d,y),x_d,y)| \\
&\qquad -|\phi_{y_d}(w,y)-\phi_{y_d}(\fx'(w,x_d,y),x_d,y)|,
\end{align*}
and therefore, 
\[
|\nabla_y(\phi(x,y)-\phi(w,y)|\ge C\max\{2^{-\ell} |x_d-w_d|,|x'-\fx'(w,x_d,y)|\}.
\]
With these estimates in place we integrate by parts in the $y$ variables, noting that for a multiindex $\alpha$
\[
|\partial_y^\alpha \sigma |\le C_{|\alpha|} 2^{\ell |\alpha|}
\]
and for $|\alpha|>1$,
\[
|\partial_y^\alpha \phi|\le C_{|\alpha|} |x-w|.
\]
Thus we get the estimate
\[
|K^{\cA\cA^*}(x,w)|\lesssim_N \int \frac{1}{(1+2^{k-\ell}|x'-\fx'(w,x_d,y)|)^N}\frac{1}{(1+2^{k-2\ell}|x_d-w_d|)^N} \, dy.
\]
Integrating in $x$ we see that
\begin{align*}
\int |K^{\cA\cA^*}(x,w)| \ dx &\lesssim_N  \int \frac{1}{(1+2^{k-2\ell}|x_d-w_d|)^N} \ dx_d \\
&\qquad \times\sup_{x_d, y} 2^{-d\ell} \int \frac{1}{(1+2^{k-\ell}|x'-\fx'(w,x_d,y)|)^N} \, dx'  \\
&\lesssim 2^{2\ell-k} 2^{-d\ell} 2^{(d-1)(\ell-k)} \\
&\lesssim 2^{\ell-dk}.
\end{align*}

\subsubsection{Proof of Lemma \ref{almostorth}}
First, (a) follows immediately since it implies 
\[
\chi(2^{\ell}y-m)\overline{\chi(2^{\ell}y-\tilde{m})}=0.
\]
The kernel $K^{\cA^*\cA}_{m,\tilde{m}}$ vanishes under the assumption in (b) because $\pi_L$ is a fold, meaning that $V_L \det d\pi_L$ is bounded away from 0 on $\cL$. Since $|\det\phi(x,y)|$ and $|\det\phi(x,z)|$ are both bounded above by $2^{-\ell+1}$, their sum is bounded by $2^{-\ell+3}$. Expanding the difference about $y=z$ we see
\begin{align*}
    \det\phi_{xy}(x,y)-\det\phi_{xy}(x,z)&=(y_d-z_d)\partial_{y_d}\det\phi_{xy}(x,z) \\
    &\qquad +(y'-z')\cdot \nabla_{y'}\det\phi_{xy}(x,z)+O(|y-z|^2) \\
    &=(y_d-z_d)[\partial_{y_d}-V_L]\det\phi_{xy}(x,z) \\
    &\qquad +(y_d-z_d)V_L\det\phi_{xy}(x,z) \\
    &\qquad +(y'-z')\cdot \nabla_{y'}\det\phi_{xy}(x,z)+O(|y-z|^2) \\
    |\det\phi_{xy}(x,y)-\det\phi_{xy}(x,z)|&\ge \frac{c_L}{3}|y_d-z_d| \\
    &\ge 5(2^{-\ell}).
\end{align*}
Thus we see there are no $y,z$ that satisfy these conditions, hence 
\[
a_{k,\ell,\pm}(x,y)\overline{a_{k,\ell,\pm}(x,z)}=0.
\]

To prove (c) we split into two cases: first, assume that $k\ge (2+\varepsilon)\ell$. Then we use the following Taylor approximation of the derivative of the phase of $K^{\cA^*\cA}_{m,\tilde{m}}$. 
\begin{equation}
    \nabla_{x'}[\phi(x,y)-\phi(x,z)]=\phi_{x'y_d}(x,z)(y_d-z_d)+\phi_{x'y'}\cdot(y'-z')+O(|y-z|^2).
\end{equation}
We know that $|\phi_{x'y'}(x,z)\cdot(y'-z')|\ge C_{d}|y'-z'|$, and $|\phi_{x'y_d}(x,z)(y_d-z_d)|\le \varepsilon|y_d-z_d|$.
By assumption 
\[
|y'-z'|\ge \tfrac{c_L}{3\|\phi\|_{C^3}}|y_d-z_d|\ge \tfrac{10\varepsilon}{C_d}|y_d-z_d|.
\]
Thus 
\begin{align*}
    |\nabla_{x'}[\phi(x,y)-\phi(x,z)]|\ge c|y-z|
\end{align*}
for some small constant $c>0$. Define the operator
\begin{equation*}
    \cM_{x'}=\tfrac{1}{i2^k}\frac{\nabla_{x'}[\phi(x,y)-\phi(x,z)]}{|\nabla_{x'}[\phi(x,y)-\phi(x,z)]|^2}\cdot \nabla_{x'}.
\end{equation*}
We apply $\cM_{x'}$ many times to $K^{\cA^*\cA}_{m,\tilde{m}}$, and by our lemma
\begin{align*}
    |K^{\cA^*\cA}_{m,\tilde{m}}(y,z)|&=\left|\int e^{i2^k(\phi(x,y)-\phi(x,z))} \sigma \enskip dx\right| \\
    &=\Big|\int e^{i2^k(\ophi(x,y)-\ophi(x,z))} \big(\cM^*_{x'}\big)^N\sigma \enskip dx\Big| \\
    &\lesssim_N \int \frac{1}{(2^{k-\ell}|y-z|)^N} |\tilde{\sigma}| \enskip dx \\
    &\lesssim_N \frac{1}{(2^{k-\ell}|y-z|)^N}\chi_m(y)\chi_{\tilde{m}}(z).
\end{align*}
Since $|y-z|>2^{-\ell}$, $k\ge 2\ell$, and $|y-z|\simeq 2^{-\ell}|m-\tilde{m}|$,
\[
\frac{1}{2^{k-\ell}|y-z|}\le \min\Big\{\frac{C}{1+2^{k-\ell}|y-z|}, \ \frac{C}{2^{k-2\ell}|m-\tilde{m}|}\Big\}
\]
Integrating in $y$ (or $z$)
\begin{align*}
\int |K^{\cA^*\cA}_{m,\tilde{m}}(y,z)| \, dy&\le C_{N,d} \int \frac{1}{(1+2^{k-\ell}|y-z|)^{d+1}} \frac{1}{(2^{k-2\ell}|m-\tilde{m}|)^N} \chi_m(y)\chi_{\tilde{m}}(z) \, dy \\
&\le C_{N,d} 2^{d(\ell-k)} (2^{k-2\ell}|m-\tilde{m}|)^{-N}.
\end{align*}
Since $k-2\ell\ge \ell\varepsilon$, if we let $N=d/\varepsilon$ then by Schur's Lemma
\[
\|\cA^*_m \cA_{\tilde{m}}\|_{2\to 2}\le C(\varepsilon,d) 2^{(\ell-dk)}|m-\tilde{m}|^{-N},
\]
proving part (c) of Lemma \ref{almostorth}.

\section{$L^p$-Sobolev Estimate}\label{PRSsection}

As in \cite{Be19}, we prove Theorem \ref{mainthm} by applying a special case of Theorem 1.1 from \cite{PrRoSe11} and a Littlewood-Paley estimate adapted from \cite{Se93}.

Let
\[
\cR_\ell=\sum_{k\ge 2\ell} \cR_{k,\ell}.
\]
We will prove for compactly supported $f$
\[
\|\cR_{\ell}f\|_{F^{p,q}_{1/p}}\le 2^{-\ell\varepsilon(p)}\|f\|_{B^{p,p}_{0}}, \qquad 0< q\le 2<4<p<\infty,
\]
where $F_{s}^{p,q}$ and $B_s^{p,q}$ are respectively the Triebel-Lizorkin space and Besov spaces (see \cite{Tr83}).
Summing in $\ell$ with $q\ge 1$ we conclude that
\[
\cR: B^{p,p}_{s,comp}\to F^{p,q}_{s+1/p}, \qquad q\le 2<4<p<\infty.
\]
Since $L^p_s=F^{p,2}_s\xhookrightarrow{} B^{p,p}_s$ for $p>2$ and $F^{p,q}_{s+1/p}\xhookrightarrow{} F^{p,2}_{s+1/p}=L^p_{s+1/p}$ for $q\le 2$, this implies the asserted $L^p$-Sobolev bounds for $\cR$. 

Let $P_k$ be standard Littlewood-Paley multipliers on $\rr^3$ for $k\in\nn$ and $\tilde{\Phi}_j(x,y)=S^j-y_j$ for $j=1,2.$ Because $\nabla_{x}\tilde{\Phi}_j(x,y)$ are linearly independent, as are 
$\nabla_{y}\tilde{\Phi}_j(x,y)$, we can find $C_0>0$ such that 
\[
4C_0^{-1}|\tau|\le |(\tau\cdot \tilde{\Phi})_x|,|(\tau\cdot\tilde{\Phi})_y|\le C_0/4 |\tau|
\]
This implies the following.
\begin{lem}\label{LittlewoodPaley} Suppose $k',k''\in \nn$, $k'\ge 2\ell$ and $\max\{|k-k'|,|k-k''|\}\ge C_1$, where $C_1$ depends on $C_0$. Then 
\[
\|P_k \cR_{k',\ell} P_{k''}\|_{L^p\to L^p}\le C \min\{2^{-kN},2^{-k'N},2^{-k''N}\}.
\]
\end{lem}

\bp[Proof of Lemma \ref{LittlewoodPaley}] This integration by parts argument is essentially due to H\"ormander \cite{Ho71}, based on the fact that the canonical relation stays away from zero sections (cf. Lemma 2.1 in \cite{Se93}). Note that the kernel of the operator $P_k \cR_{k',\ell} P_{k''}$ is given by
\begin{align*}
\int\int\int\int\int e^{i\left[\langle x-w,\eta\rangle+\tau\cdot\tilde{\Phi}(w,z)+\langle z-y,\xi\rangle\right]} \chi_0(2^{-k}|\eta|)\chi_0(2^{-k'}|\tau|)\chi_0(2^{-k''}|\xi|) \\
\qquad \times a_{k,\ell,\pm}(z_1,\tau) \chi(|w|)\chi(|x|)\, dw \,dz \,d\tau \,d\eta \,d\xi.
\end{align*}
Our assumption on $\Phi$ implies that if $\max\{|k-k'|,|k'-k''|\}>C_1$ we have
\[
\nabla_{(z,w)}\left[ \langle x-w,\eta\rangle+\tau\cdot\tilde{\Phi}(w,z)+\langle z-y,\xi\rangle\right]\ge c \max\{2^k,2^{k'},2^{k''}\}.
\]
Thus we integrate by parts in the $(w,z)$ variables to get the above bound on the kernel, implying by Minkowski the desired bound on $L^p$. 

\ep

Using the lemma above and an argument similar to a part of the proof of Lemma 2.1 in \cite{Se93}, we can reduce the proof of Theorem \ref{mainthm} to the estimate
\begin{equation}\label{PRSinequality}
    \Big\|\Big(\sum_{k\ge 2\ell} \big|2^{k/p}P_k\cR_{k+s_1,\ell}P_{k+s_2} f\big|^q\Big)^{1/q}\Big\|_{L^p}\le C 2^{-\ell\varepsilon(p)}\Big\|\Big(\sum_{k>0} |P_{k+s_2}f|^p\Big)^{1/p}\Big\|_{L^p}.
\end{equation}

To prove \eqref{PRSinequality} we apply the main result from \cite{PrRoSe11}.

\begin{thm}[\cite{PrRoSe11}]\label{PRS}
Let $T_k$ be a family of operators defined for Schwartz functions by
\[
    T_k f(x)=\int K_k(x,y) f(y) \, dy.
\]
Let $\phi\in \mathcal{S}(\rr^3)$, $\phi_k=2^{3k}\phi(2^k\cdot)$, and $\Pi_k f=\phi_k*f$. Let $\varepsilon>0$ and $1<p_0<p<\infty$. Assume $T_k$ satisfies
\begin{align}\label{PRSLp}
    \sup_{k>0} 2^{k/p}\|T_k\|_{L^p\to L^p}&\le A \\
    \sup_{k>0} 2^{k/p_0}\|T_k\|_{L^{p_0}\to L^{p_0}}&\le B_0. \label{PRSLpo}
\end{align}
Further let $A_0 \ge 1$, and assume that for each cube $Q$ there is a measurable set $E_Q$ such that
\begin{equation}\label{EQbound}
    |E_Q|\le A_0\max\{|Q|^{2/3},|Q|\},  
\end{equation}
and for every $k\in\nn$ and every cube $Q$ with $2^k {\rm diam }(Q)\ge 1,$
\begin{equation}\label{EQintegral}
    \sup_{x\in Q} \int_{\rr^d\setminus E_Q} |K_k(x,y)| \, dy\le B_1 \max\left\{\left(2^k{\rm diam}(Q)\right)^{-\varepsilon},2^{-k\varepsilon}\right\}.
\end{equation}
Let 
\[
\mathcal{B}=B_0^{q/p}(A A_0^{1/p}+B_1)^{1-q/p}.
\]
Then for any $q>0$ there is a $C$ depending on $\varepsilon,p,p_0,q$ such that
\begin{equation}\label{PRSestimate}
    \Big\|\Big(\sum_k 2^{kq/p}|P_k T_k f_k|^q\Big)^{1/q}\Big\|_{p}\le C A \left[\log\left(3+\tfrac{\mathcal{B}}{ A}\right)\right]^{1/q-1/p}\Big(\sum_k \|f_k\|_p^p\Big)^{1/p}.
\end{equation}
\end{thm}

We apply this theorem on the family of operators $T_k:=\cR_{k,\ell}$ for $k\ge 2\ell$ (here $\ell$ is fixed). By Proposition \ref{decomposedfinal} the assumptions \eqref{PRSLp} and \eqref{PRSLpo} are satisfied with $A\lesssim 2^{-\ell\varepsilon(p)}$ and $B_0\lesssim 2^{-\ell\varepsilon(p_0)}$. We next check assumptions \eqref{EQbound} and \eqref{EQintegral}. For a given cube $Q$ with center $x_Q$ let
\[
E_Q=\{y \, : \, |S(x^Q,y_3)-y'|\le C 2^{\ell}\diam(Q)\}
\]
if $\diam(Q)<1$, and a cube centered at $x^Q$ of diameter $C2^{\ell}\diam(Q)$ if $|Q|\ge 1$. By an integration by parts argument we derive the bound
\[
|K_k(x,y)|\lesssim_N\frac{2^{2k}}{(1+2^{k-\ell}|S(x^Q,y_3)-y'|)^N}.
\]
Then clearly assumptions \eqref{EQbound} and \eqref{EQintegral} are satisfied with $A_0\lesssim 2^{3\ell}$ and $B_1\lesssim 2^{2\ell}$ respectively. Theorem \ref{PRS} then implies \eqref{PRSinequality} with $\Pi_k=P_{k+s_1}$ and $f_k=P_{k+s_2} f$, finishing the proof of Theorem \ref{mainthm}.

\bibliography{mybib}
	\bibliographystyle{plain}

\end{document}